\documentclass[a4paper,10pt]{amsart}

\usepackage[english]{babel}
\usepackage[latin1]{inputenc}
\usepackage{amsmath}
\usepackage{amsthm}
\usepackage{amsfonts}
\usepackage{amssymb, latexsym, graphics, showidx}

\numberwithin{equation}{section}

\newtheorem{thm}{Theorem}[section]
\newtheorem{rem}[thm]{Remark}
\newtheorem{cor}[thm]{Corollary}
\newtheorem{prop}[thm]{Proposition}
\newtheorem{lem}[thm]{Lemma}
\newtheorem*{ack}{Acknowledgments}

\renewcommand{\dim}{\noindent\textbf{Proof.} }
\newcommand{\dims}{\noindent\textbf{Proof} }
\newcommand{\finedim}{{\unskip\nobreak\hfil\penalty50
   \hskip2em\hbox{}\nobreak\hfil\mbox{$\Box$ \qquad}
   \parfillskip=0pt \finalhyphendemerits=0\par\medskip}}
\newcommand{\R}{\mathbb{R}}
\newcommand{\N}{\mathbb{N}}
\newcommand{\Om}{\Omega}
\newcommand{\Oms}{\overline{\Omega}}
\newcommand{\lam}{\lambda}
\newcommand{\al}{\alpha}
\newcommand{\lams}{\overline{\lambda}}
\newcommand{\lamso}{\underline{\lambda}}
\newcommand{\xs}{\overline{x}}
\newcommand{\ys}{\overline{y}}
\newcommand{\zs}{\overline{z}}

\title[The Neumann eigenvalue problem]
{Principal eigenvalues for Isaacs operators with Neumann boundary
conditions}
\author{Stefania Patrizi}
\address{SAPIENZA Universit\`a di Roma, Dipartimento di Matematica,
Piazzale A.~Moro 2, I-00185 Roma, Italy}
\email{patrizi@mat.uniroma1.it}

\begin{document}
\keywords{Fully nonlinear uniformly elliptic operators, Isaacs
operators, Neumann boundary value problems, Maximum Principle,
principal eigenvalues, viscosity solutions.}

\begin{abstract}
In this paper we show the existence of two principal
eigen\-va\-lues associated to general non-convex fully nonlinear
elliptic operators with Neumann boun\-da\-ry conditions in a
bounded $C^2$ domain. We study these objects and we establish some
of their basic properties. Finally, Lipschitz regularity,
uniqueness and existence results for the solution of the Neumann
problem are given.
\end{abstract}

\maketitle

\section{Introduction}

A self-adjoint uniformly elliptic linear operator in a bounded
domain possesses a countable set of real eigenvalues, the least of
them, called first or principal eigenvalue, can be
cha\-rac\-te\-rized as the minimum of the associated Rayleigh
quotient, see e.g. \cite{gt}. This characterization is possible
only for operators which are variational. In their famous work
\cite{bnv}, Berestycki, Nirenberg and Varadhan defined the
principal eigenvalue $\lam_1$ of a general uniformly elliptic
linear operator
$$L[u]=-\text{tr}(A(x)D^2u)+b(x)\cdot Du+c(x)u,$$
in a bounded domain $\Om$, as the supremum of those $\lam$ for
which there exists a positive supersolution of $L[u]=\lam u$. In
that paper, they showed that $\lam_1$ is the least eigenvalue of
$L$, i.e., for any eigenvalue $\lam\neq \lam_1$,
Re$(\lam)>\lam_1$; moreover $\lam_1$ can be characterized as the
supremum of those $\lam$ for which the operator $L-\lam I$
satisfies the maximum principle, i.e., for any $\lam<\lam_1$, if
$u$ is a subsolution of $L[u]-\lam u=0$ and $u\leq 0$ on
$\partial\Om$ then $u\leq 0$ in $\Om$. Furthermore, they
established several other properties of the first eigenvalue, such
as simplicity and stability.

In view of its relation with the maximum and the comparison
principles, the concept of principal eigenvalue has been extended
to nonlinear operators to study the associated boundary value
problems. That has been done for the variational operators, such
as the p-Laplacian, through the method of minimization of the so
called nonlinear Rayleigh quotient,  see e.g. \cite{an} and
\cite{li2}. An important step in the study of the eigenvalue
problem for general nonlinear operators was made by Lions in
\cite{l}. In that paper, using probabilistic and analytical
methods, he showed the existence of principal eigenvalues for the
uniformly elliptic Hamilton-Jacobi-Bellman operator
\begin{equation}\label{hjb}F(x,u,Du,D^2u)=\sup_{\al\in \mathcal{A}}\{ -\text{tr}(A_\al(x) D^2u)+b_\al(x)\cdot Du+c_\al(x)
u\},\end{equation} which arises in stochastic optimal control, see
\cite{bc} and \cite{fs}.  Very recently, many authors, inspired by
\cite{bnv}, have de\-ve\-lo\-ped an eigenvalue theory for fully
nonlinear operators which are non-variational. The Pucci's
extremal operators $\mathcal{M}_{a,A}(D^2u)$ (see the next
section) have been treated by  Quaas \cite{q} and Busca, Esteban
and Quaas \cite{beq}. Related results have been obtained by
Birindelli and Demengel in \cite{bd} for singular or degenerate
elliptic operators, like $|Du|^\al\mathcal{M}_{a,A}(D^2u)$,
$\al>-1$, the p-Laplacian and some its non-variational
generalizations. In \cite{qs} Quaas and Sirakov have studied the
eigenvalue problem for fully nonlinear elliptic operators which
are convex and positively homogenous, like the
Hamilton-Jacobi-Bellman one \eqref{hjb}. In that paper many
properties of the principal eigenvalues, including the fact that
they are simple and isolated, have been established. Similar
results have been obtained by Ishii and Yoshimura \cite{iy} for
non-convex operators, such as the Isaacs one
\begin{equation}\label{isaacs}F(x,u,Du,D^2u)=\sup_{\al\in \mathcal{A}}\inf_{\beta\in \mathcal{B}}
\{- \text{tr}(A_{\al,\beta}(x) D^2u)+b_{\al,\beta}(x)\cdot
Du+c_{\al,\beta}(x) u\},\end{equation}which arises in stochastic
differential games, see \cite{fs}.

All these articles treat Dirichlet boundary conditions.

In this paper we want to develop an eigenvalue theory for a class
of fully nonlinear operators with Neumann boundary conditions in a
bounded $C^2$ domain $\Om$. Precisely, we consider a uniformly
elliptic operator which is positively homogenous of order 1
\begin{equation}\label{e1}F[u](x)=F(x,u,Du,D^2u),
 \end{equation}for any $u\in C^2(\Oms)$, with some additional
 assumptions that will be made precise in the next section. This
 class includes the non-convex Isaacs operator \eqref{isaacs} if
 $c_{\al,\beta}$ are equi-continuous and
 $A_{\al,\beta}$, $b_{\al,\beta}$  are
 equi-Lipschitz continuous.

 To \eqref{e1} we associate the following boundary condition
 \begin{equation}\label{e2}
B(x,u,Du)=f(x,u)+\langle Du,\overrightarrow{n}(x)\rangle=0\quad
x\in\partial\Om,
\end{equation}where $\overrightarrow{n}(x)$ is the exterior normal to the
domain $\Om$ at $x$. Typically $f$ will be
\begin{equation}\label{flinear}f(x,u)=\gamma(x)u,\end{equation} with $\gamma(x)$ continuous
 and non-negative on $\partial\Om$.

Following the ideas of \cite{bnv}, we define the principal
eigenvalues as
\begin{equation*}\begin{split}\lams:=\sup\{&\lam\in
\R\;|\;\exists\,v>0 \text{ on $\Oms$ bounded viscosity
supersolution of }\\&
F(x,v,Dv,D^2v)=\lam v \text{ in } \Om,\, B(x,v,Dv) = 0\text{ on }\partial \Om \},\\
\end{split}\end{equation*}
\begin{equation*}\begin{split}\lamso:=\sup\{&\lam\in
\R\;|\;\exists\,u<0 \text{ on $\Oms$ bounded  viscosity
subsolution of }\\&
F(x,u,Du,D^2u)=\lam u \text{ in } \Om,\, B(x,u,Du)= 0\text{ on }\partial \Om \}.\\
\end{split}\end{equation*} One of the
scope of this work is to prove that $\lams$ and $\lamso$ are
"eigenvalues" for $F$ which admit respectively a positive and a
negative "eigenfunction".
 Moreover, we  show that $\lams$ (resp., $\lamso$) can be characterized as the supremum of those
$\lam$ for which the ope\-ra\-tor $F-\lam I$ with boundary
condition \eqref{e2} satisfies the maximum (resp., minimum)
principle. As a consequence, $\lams$ (resp., $\lamso$) is the
least "eigenvalue" to which there cor\-re\-spond "eigenfunctions"
positive (resp., negative) somewhere.

Other properties of the principal eigenvalues are established: we
 show that they are simple, isolated and the only
"eigenvalues" to which there correspond "eigenfunctions" which do
not change sign in $\Om$. Finally, we obtain Lipschitz regularity,
uniqueness and existence results for viscosity solutions of
\begin{equation}\label{sisintro}
\begin{cases}
F(x,u,Du,D^2u)=g(x) & \text{in} \quad\Om \\
 B(x,u,Du)= 0 & \text{on} \quad\partial\Om. \\
 \end{cases}
 \end{equation}
In particular, we prove that \eqref{sisintro} is solvable for any
continuous right-hand side if the two principal eigenvalues are
positive.

The paper is organized as follows. In the next section we give
assumptions and precise the notion of solution adopted. In Section
3 we prove the strong comparison principle between sub and
supersolutions of \eqref{sisintro}. This allows us to prove the
maximum principle for subsolutions of the Neumann boundary value
problem. We first show it under the classical assumption that $F$
be proper, see Theorem \ref{macpcless0}; then we prove in Theorem
\ref{maxpneum} that the operator $F-\lam I$ with boundary
condition \eqref{e2} satisfies the maximum principle for any
$\lam<\lams$. Using the example given in \cite{p} we show that the
result of Theorem \ref{maxpneum} is stronger than that of Theorem
\ref{macpcless0}, i.e., that there exist non-proper operators
which have positive principal eigenvalue $\lams$, and then for
which the maximum principle holds.

In Section 5 we establish a Lipschitz regularity result for
viscosity solutions of \eqref{sisintro}. In Section 6 we show some
existence and comparison theorems. In Section 7 we establish some
of the basic properties of the principal
 eigenvalues. Finally, in Section 8 we show, through an example, that
 $\lams$ and $\lamso$ may be different.

In \cite{p} the author of this paper has studied the principal
eigenvalues of fully nonlinear singular elliptic operators modeled
on the p-Laplacian or on $|Du|^\al\mathcal{M}_{a,A}(D^2u)$,
$\al>-1$, with the pure Neumann boundary condition. In that paper
we have used a different approach since the strong comparison
principle between sub and supersolutions is not known. For this
reason some questions about the properties of the principal
eigenvalues, such as simplicity, were left open.
\section{Assumptions}
Let $\emph{S(N)}$ be the space of symmetric matrices on $\R^N$,
equipped with the usual ordering. We denote by
$\mathcal{M}_{a,A}^+$, $\mathcal{M}_{a,A}^-:\emph{S(N)}\rightarrow
\R$ the Pucci's extremal operators defined by
$$\mathcal{M}_{a,A}^+(X)=A\sum_{e_i>0}e_i+a\sum_{e_i<0}e_i,$$
$$\mathcal{M}_{a,A}^-(X)=a\sum_{e_i>0}e_i+A\sum_{e_i<0}e_i,$$where
$e_1,...,e_N$ are the eigenvalues of $X$ (see e.g. \cite{cc}).

The operator $F$ is supposed to be continuous on
$\overline{\Om}\times\R\times\R^N \times \emph{S(N)}$, moreover we
shall make the following assumptions:
\renewcommand{\labelenumi}{(F\arabic{enumi})}
\begin{enumerate}
 \item For all $(x,r,p,X)\in\overline{\Om}\times\R\times\R^N
\times \emph{S(N)}$ and $t \geq
 0$ $$F(x,tr,tp,t X)=tF(x,r,p,X).$$
 \item There exist $b,c>0$ such that for $x\in\overline{\Om},\,r,s\in\R,\, p,q\in \R^N,\, X, Y\in \emph{S(N)}$
 \begin{equation*}\begin{split} \mathcal{M}_{a,A}^-(Y-X)-b|p-q|-c|r-s|&\leq
 F(x,r,p,X)-F(x,s,q,Y)\\&
 \leq \mathcal{M}_{a,A}^+(Y-X)+ b|p-q|+c|r-s|.\end{split}\end{equation*}
 \item For each $T>0$ there exists a continuous function $\omega_T$ with $\omega_T(0)=0$, such that if  $X,Y\in \emph{S(N)}$ and $\zeta>0$ satisfy\\
 \begin{equation*}-3\zeta\left(%
\begin{array}{cc}
  I & 0 \\
  0 & I \\
\end{array}%
\right)\leq \left(%
\begin{array}{cc}
  X & 0 \\
  0 & -Y \\
\end{array}%
\right)\leq 3\zeta\left(%
\begin{array}{cc}
  I & -I \\
  -I & I \\
\end{array}%
\right)
\end{equation*}
and $I$ is the identity matrix in $\R^N$, then for all
$x,y\in\Oms$, $r\in [-T,T]$, $p\in\R^N$
$$F(y,r,p,Y)-F(x,r,p,X)\leq \omega_T (\zeta|x-y|^2+|x-y|(|p|+1)).$$
  \item There exists $C_1>0$ such that for all $x,y\in\Oms$ and $X\in\emph{S(N)}$
 $$|F(x,0,0,X)-F(y,0,0,X)|\leq C_1|x-y|^\frac{1}{2}\|X\|.$$
\end{enumerate}
Here and in what follows we fix the norm $\|X\|$ in $\emph{S(N)}$
by setting
$$\|X\|=\sup\{|X\xi|\,:\,\xi\in\R^N,\,|\xi|\leq1\}=\sup\{|\lam|\,:\,
\lam \text{ is an eigenvalue of }X\}.$$ Remark that (F1) implies
that $F(x,0,0,0)\equiv 0$.

The Isaacs operator \eqref{isaacs} is continuous and satisfies
(F1) and (F2) if $aI\leq A_{\al,\beta}(x)\leq AI$ for any
$x\in\Oms$, $(\al,\beta)\in \mathcal{A}\times\mathcal{B}$ and the
functions $A_{\al,\beta},\,b_{\al,\beta},\,c_{\al,\beta}$ are
continuous on $\Oms$ uniformly in $\al$ and $\beta$, where
$\mathcal{A}$ and $\mathcal{B}$ are arbitrary index sets. If the
matrices $A_{\al,\beta}$ are equi-H\"{o}lderian of exponent
$\frac{1}{2}$, i.e., for some constant $C>0$
$$\|A_{\al,\beta}(x)-A_{\al,\beta}(y)\|\leq C|x-y|^\frac{1}{2}\quad \text{for all }x,y\in\Oms\text{ and }(\al,\beta)\in\mathcal{A}\times\mathcal{B},$$
then $F$ satisfies (F4). Finally, (F3) is satisfied by $F$ if, in
addition to the uniform elliptic condition $A_{\al,\beta}(x)\geq
aI$ and the equi-continuity of $c_{\al,\beta}$, the functions
$A_{\al,\beta}$ and $b_{\al,\beta}$ are equi-Lipschitz continuous,
i.e., there exists $L>0$ such that for all $x,y\in\Oms$ and
$(\al,\beta)\in\mathcal{A}\times\mathcal{B}$
$$\|A_{\al,\beta}(x)-A_{\al,\beta}(y)\|\leq L|x-y|,\quad|b_{\al,\beta}(x)-b_{\al,\beta}(y)|\leq L|x-y|.$$

We assume throughout the paper that $\Om$ is a bounded domain of
$\R^N$ of class $C^2$. In particular it satisfies the interior
sphere condition and the uniform exterior sphere condition, i.e.,
\begin{itemize}
\item [($\Om1$)]For
 each $x\in\partial \Om$ there exist $R>0$ and $y\in \Om$ for
 which $|x-y|=R$ and $B(y,R)\subset \Om$.
\item[($\Om2$)] There exists $r>0$ such that  $B(x+r\overrightarrow{n}(x),r)\cap \Om =\emptyset$  for any $x\in \partial\Om.$
\end{itemize}
From property ($\Om$2) it follows that
\begin{equation}\label{sferaest}
\langle
\overrightarrow{n}(x),y-x\rangle\leq\frac{1}{2r}|y-x|^2\quad\text{for
}x\in\partial\Om \text{ and } y\in\Oms.\end{equation} Moreover,
the $C^2$-regularity of $\Om$ implies the existence of a
neighborhood of $\partial\Om$ in $\Oms$ on which the distance from
the boundary
 $$d(x):=\inf\{|x-y|, y\in\partial\Om\},\quad x\in\Oms$$ is of class $C^2$. We still denote by $d$ a
$C^2$ extension of the distance function to the whole $\Oms$.
Without loss of generality we can assume that $|Dd(x)|\leq1$ on
$\Oms$.

On the function $f$ in \eqref{e2} we shall suppose
\begin{itemize}
\item [(f1)] $f:\partial\Om\times\R\rightarrow\R$ is continuous.
\end{itemize}
\begin{itemize}
\item [(f2)]For all $(x,r)\in\partial\Om\times\R$ and $t\geq 0$ $$f(x,tr)=tf(x,r).$$
\end{itemize}
For the existence results we will assume in addition
\begin{itemize}
\item [(f3)]For all $x\in\partial\Om$ $r\rightarrow f(x,r)$ is non-decreasing on $\R$.
\end{itemize}
Clearly,  $f(x,u)=\gamma(x)u$  with $\gamma(x)\geq0$ and
continuous on $\partial\Om$, satisfies all the three hypothesis.

In this paper we adopt the notion of viscosity solution. We denote
by $USC(\overline{\Om})$ the set of upper semicontinuous functions
on $\overline{\Om}$ and by $LSC(\overline{\Om})$ the set of lower
semicontinuous functions on $\overline{\Om}$. Given
$g:\Oms\rightarrow\R$, we recall that a function $u\in
USC(\overline{\Om})$ (resp., $u\in LSC(\overline{\Om})$) is called
\emph{viscosity subsolution} (resp., \emph{supersolution}) of
\begin{equation}\label{syscap2}
\begin{cases}
 F(x,u,Du,D^2u)= g(x) & \text{in} \quad\Om \\
 B(x,u,Du)= 0 & \text{on} \quad\partial\Om. \\
 \end{cases}
 \end{equation}if the following conditions hold
\begin{itemize}
\item[(i)] For every $x_0\in \Om$ for all $\varphi\in C^2(\overline{\Om})$, such that $u-\varphi$ has a local maximum (resp., minimum)
on $x_0$ then
$$F(x_0,u(x_0),D\varphi(x_0),D^2\varphi(x_0))\leq \,(\text{resp., } \geq\,)\,g(x_0).$$
\item[(ii)]For every $x_0\in \partial\Om$ for all $\varphi\in C^2(\overline{\Om})$, such that $u-\varphi$ has a local maximum (resp., minimum)
on $x_0$ then
$$(F(x_0,u(x_0),D\varphi(x_0),D^2\varphi(x_0))-g(x_0))\wedge
B(x_0,u(x_0),D\varphi(x_0))\leq 0$$(resp.,
$$(F(x_0,u(x_0),D\varphi(x_0),D^2\varphi(x_0))-g(x_0))\vee B(x_0,u(x_0),D\varphi(x_0))\geq 0.)$$
\end{itemize}
A \emph{viscosity solution} is a continuous function which is both
a subsolution and a supersolution.

In the above definition the test functions can be substituted by
the elements of the semi-jets $\overline{J}^{2,+}u(x_0)$ when $u$
is a subsolution and $ \overline{J}^{2,-}u(x_0)$ when $u$ is a
supersolution. For a detailed presentation of the theory of
viscosity solutions we refer the reader to e.g. \cite{cil}.

One of the motivation for these relaxed boundary conditions is the
stability under uniform convergence. Actually, if the operator $F$
satisfies (F2) and the domain $\Om$ the exterior sphere condition,
viscosity subsolutions (resp., supersolutions) satisfy in the
viscosity sense $B(x,u(x),Du(x))\leq $ (resp. $\geq$ ) 0 for any
$x\in\partial\Om$, as shown in the following proposition due to
Hitoshi Ishii, \cite{icom}, whose proof is given for the reader's
convenience.

\begin{prop}\label{ishii}Suppose that $\Om$ satisfies the exterior sphere
condition. If there exists $b>0$ such that for
$x\in\overline{\Om},\,r\in\R,\, p,q\in \R^N,\, X, Y\in
\emph{S(N)}$ $$
F(x,r,p,X)-F(x,r,q,Y)\geq\mathcal{M}_{a,A}^-(Y-X)-b|p-q|,$$ and
$u$ is a viscosity subsolution of \eqref{syscap2} then  $u$
satisfies in the viscosity sense
$$B(x_0,u(x_0),Du(x_0))\leq 0,$$ for any
$x_0\in\partial \Om$. If
$$F(x,r,p,X)-F(x,r,q,Y)\leq\mathcal{M}_{a,A}^+(Y-X)+b|p-q|,$$ and $u$ is a viscosity supersolution of \eqref{syscap2} then  $u$
satisfies in the viscosity sense
$$B(x_0,u(x_0),Du(x_0))\geq 0,$$ for any
$x_0\in\partial \Om$.
\end{prop}
\dim We show the proposition for subsolutions. Set
$$g(t)=-Kt^2+\epsilon t\quad \forall t\in\R,$$ where $K\gg1$ and
$0<\epsilon\ll1$. Observe that $g(0)=0$, $g'(0)=\epsilon$,
$g''(0)=-2K$, and
$$0<t<\frac{\epsilon}{K}\quad\Longrightarrow\quad g(t)>0.$$
Let $\varphi\in C^2(\Oms)$ and $x_0\in\partial\Om$. Assume that
$u-\varphi$ attains a maximum at $x_0$. We need to prove that
$f(x_0,u(x_0))+\langle
\overrightarrow{n}(x_0),D\varphi(x_0)\rangle\leq0$.

Let $y_0\in\R^N$ and $R>0$ satisfy
$$B(y_0,R)\cap\Oms=\{x_0\}.$$ We may assume by translation that
$y_0=0$. We set
$$\psi(x)=g(|x|-R)\quad \forall x\in\R^N.$$ Note that
$\psi(x_0)=g(0)=0,$
$$D\psi(x_0)=g'(0)\frac{x_0}{|x_0|}=\epsilon e_0,\text{ where
}e_0=\frac{x_0}{|x_0|},$$ $$\overrightarrow{n}(x_0)\cdot
D\psi(x_0)=-e_0\cdot\epsilon e_0=-\epsilon,$$
\begin{equation*}\begin{split}D^2\psi(x_0)&=g''(0)e_0\otimes
e_0+\frac{g'(0)}{|x_0|}(I-e_0\otimes e_0)\\&=-2Ke_0\otimes
e_0+\frac{\epsilon}{R}(I-e_0\otimes
e_0),\end{split}\end{equation*}
$$\mathcal{M}_{a,A}^-(-D^2\psi(x_0))=-\frac{\epsilon(N-1)A}{R}+2Ka,$$
$$R<|x|<R+\frac{\epsilon}{K}\quad\Longrightarrow\quad \psi(x)>0.$$
Moreover we observe that $u-\varphi-\psi$ attains a local maximum
at $x_0$. Remark that
\begin{equation*}\begin{split}
&F(x_0,u(x_0),D\varphi(x_0)+D\psi(x_0),D^2\varphi(x_0)+D^2\psi(x_0))\\&
\geq
F(x_0,u(x_0),D\varphi(x_0),D^2\varphi(x_0))-b|D\psi(x_0)|+\mathcal{M}_{a,A}^-(-D^2\psi(x_0))\\&\geq
F(x_0,u(x_0),D\varphi(x_0),D^2\varphi(x_0))-b\epsilon-\frac{\epsilon(N-1)A}{R}+2Ka.\end{split}\end{equation*}
We fix $K\gg1$ so that for any $0<\epsilon<1$
$$F(x_0,u(x_0),D\varphi(x_0),D^2\varphi(x_0))-b\epsilon-\frac{\epsilon(N-1)A}{R}+2Ka>g(x_0).$$
Then, by definition of subsolution we get that
$$0\geq f(x_0,u(x_0))+\overrightarrow{n}(x_0)\cdot
(D\varphi(x_0)+D\psi(x_0))=f(x_0,u(x_0))+\overrightarrow{n}(x_0)\cdot
D\varphi(x_0)-\epsilon,$$ from which we obtain
$$f(x_0,u(x_0))+\overrightarrow{n}(x_0)\cdot D\varphi(x_0)\leq
0,$$as desired. \finedim

\section{The Strong Comparison Principle}
The strong comparison principle is the key ingredient in the
development of our theory.
\begin{thm}\label{stcompneu} Assume
that (F2), (F3), (f1) hold and that $g$ is continuous on $\Oms$.
Let $u\in USC(\Oms)$ and $v\in LSC(\Oms)$ be respectively a sub
and a supersolution of
\begin{equation*}
\begin{cases}
 F(x,u,Du,D^2u)= g(x) & \text{in} \quad\Om \\
 B(x,u,Du)= 0 & \text{on} \quad\partial\Om. \\
 \end{cases}
 \end{equation*}
If $u\leq v$ on $\Oms$ then either $u<v$ on $\Oms$ or $u\equiv v$
on $\Oms$.
\end{thm}
Let us recall that for the sub and the supersolutions of the
Dirichlet problem the following theorem holds, see \cite{iy}.
\begin{thm}\label{puccistcomdir}
Assume that (F2), (F3) hold and that $g$ is continuous on $\Om$.
Let $u\in USC(\Om)$ and $v\in LSC(\Om)$ be respectively a sub and
a supersolution of
\begin{equation*}
 F(x,u,Du,D^2u)= g(x).
 \end{equation*}
If $u\leq v$ in $\Om$ then either $u<v$ in $\Om$ or $u\equiv v$ in
$\Om$.
\end{thm}

\dims {\bf of Theorem \ref{stcompneu}.} Assume $u\not\equiv v$,
then by Theorem \ref{puccistcomdir} $u<v$ in $\Om$. Suppose by
contradiction that there exists a point $x_0\in\partial\Om$ on
which $u(x_0)=v(x_0)$.

The interior sphere condition ($\Om$1) implies that there exist
$R>0$ and $y_0\in\Om$ such that the ball centered in $y_0$ and of
radius $R$, $B_1$, is contained in $\Om$ and $x_0 \in\partial
B_1$. Let for $k>2/R^2$ and $x\in\Oms$
$$w(x):=e^{-kR^2}-e^{-k|x-y_0|^2}.$$ This function has the
following properties \begin{equation*}\begin{split}&w(x)<0\quad
\text{in }B_1,\\
&w(x)=0\quad \text{on }\partial B_1,\\
&w(x)>0\quad\text{outside
}\overline{B}_1.\end{split}\end{equation*}Let $B_2$ be the ball of
center $y_0$ and radius $\frac{R}{2}$ and
$-m:=\max_{\overline{B}_2}(u-v)<0$. Choose $\sigma>0$ so small
that
\begin{equation}\label{sigmaw}\sigma\inf_{\overline{B}_2}
w\geq -\frac{m}{2}.\end{equation} Let us define for $j\in\N$ the
functions

$$\phi(x,y):=\frac{j}{2}|x-y|^2+\frac{\sigma}{2}(w(x)+w(y))-f(x_0,u(x_0))\langle \overrightarrow{n}(x_0),x-y\rangle,$$and
$$\psi(x,y):=u(x)-v(y)-\phi(x,y).$$
Let $(x_j,y_j)\in\Oms^2$ be a maximum point of $\psi$ in $\Oms^2$.
We have
\begin{equation}\label{stcompuv}\begin{split}0=u(x_0)-v(x_0)-\sigma w(x_0)&\leq
u(x_j)-v(y_j)-\frac{j}{2}|x_j-y_j|^2-\frac{\sigma}{2}(w(x_j)+w(y_j))\\&+f(x_0,u(x_0))\langle
\overrightarrow{n}(x_0),x_j-y_j\rangle,\end{split}\end{equation}
from which we can see that $|x_j-y_j|\rightarrow 0$ as
$j\rightarrow+\infty$. Up to subsequence, $x_j$ and $y_j$ converge
to some $\zs\in\Oms$. Standard arguments show that
$$\lim_{j\rightarrow+\infty}\frac{j}{2}|x_j-y_j|^2=0,\quad
 \lim_{j\rightarrow+\infty}u(x_j)\rightarrow u(\zs)\text{ and } \lim_{j\rightarrow+\infty}v(y_j)\rightarrow v(\zs).$$ Passing to
the limit in \eqref{stcompuv} we get
\begin{equation}\label{stcompstimau-v}\sigma w(\zs)\leq
u(\zs)-v(\zs)\leq0,\end{equation}which implies that the limit
point $\zs$ belongs to $\overline{B}_1$. Furthermore, since
$u(\zs)-v(\zs)-\sigma w(\zs)\geq 0$, it cannot belong to
$\overline{B}_2$, indeed by \eqref{sigmaw} we have
$u(x)-v(x)-\sigma w(x)\leq -\frac{m}{2}<0$, for any $x\in
\overline{B}_2$. In conclusion $$\frac{R}{2}<|\zs-y_0|\leq R.$$

Computing the derivatives of $\phi$ we get
$$D_x\phi(x,y)=j(x-y)+\sigma ke^{-k|x-y_0|^2}(x-y_0)-f(x_0,u(x_0))
\overrightarrow{n}(x_0),$$
$$D_y\phi(x,y)=-j(x-y)+\sigma ke^{-k|y-y_0|^2}(y-y_0)+f(x_0,u(x_0))
\overrightarrow{n}(x_0).$$ If $x_j\in\partial\Om$ then $\zs=x_0$
and using \eqref{sferaest} we have
\begin{equation*}\begin{split}B(x_j,u(x_j),D_x\phi(x_j,y_j)) &\geq
f(x_j,u(x_j))-f(x_0,u(x_0))
\langle\overrightarrow{n}(x_0),\overrightarrow{n}(x_j)\rangle\\&-\frac{j}{2r}|x_j-y_j|^2+\sigma
k e^{-k|x_j-y_0|^2}\langle
x_j-y_0,\overrightarrow{n}(x_j)\rangle>0\end{split}\end{equation*}
for large $j$, since the last term goes to $\sigma ke^{-kR^2}R$ as
$j\rightarrow+\infty$, being
$\overrightarrow{n}(x_0)=\frac{x_0-y_0}{R}.$ Similarly if
$y_j\in\partial\Om$ then $\zs=x_0$ and $u(x_0)=v(x_0)$ so that
\begin{equation*}\begin{split}B(y_j,v(y_j),-D_y\phi(x_j,y_j))&\leq f(y_j,v(y_j))-f(x_0,u(x_0))
\langle\overrightarrow{n}(x_0),\overrightarrow{n}(y_j)\rangle
\\&+\frac{j}{2r}|x_j-y_j|^2-\sigma k e^{-k|y_j-y_0|^2}\langle
y_j-y_0,\overrightarrow{n}(y_j)\rangle<0\end{split}\end{equation*}
for large $j$. Then $x_j$ and $y_j$ are internal points and
$$F(x_j,u(x_j),D_x\phi(x_j,y_j),X)\leq g(x_j)\quad \text{if }
(D_x\phi(x_j,y_j),X)\in \overline{J}^{2,+}u(x_j),$$
$$F(y_j,v(y_j),-D_y\phi(x_j,y_j),Y)\geq g(y_j)\quad \text{if }
(-D_y\phi(x_j,y_j),Y)\in \overline{J}^{2,-}v(y_j).$$ Then the
previous relations hold for $(x_j,y_j)\in\Oms^2$, provided $j$ is
large.

Since $(x_j,y_j)$ is a local maximum point of
$\psi(x,y)=(u(x)-\frac{\sigma}{2} w(x))-(v(y)+\frac{\sigma}{2}
w(y))-\frac{j}{2}|x-y|^2+f(x_0,u(x_0))\langle
\overrightarrow{n}(x_0),x-y\rangle$ in $\Oms^2$, applying Theorem
3.2 of \cite{cil} there exist $X_j,Y_j\in$\emph{S(N)} such that
$(D_x\phi(x_j,y_j),X_j)\in \overline{J}^{2,+}u(x_j)$,
$(-D_y\phi(x_j,y_j),Y_j)\in \overline{J}^{2,-}v(y_j)$ and
 \begin{equation*}-3j\left(%
\begin{array}{cc}
  I & 0 \\
  0 & I \\
\end{array}%
\right)\leq \left(%
\begin{array}{cc}
  X_j-\frac{\sigma}{2}D^2w(x_j) & 0 \\
  0 & -(Y_j+\frac{\sigma}{2} D^2w(y_j)) \\
\end{array}%
\right)\leq 3j\left(%
\begin{array}{cc}
  I & -I \\
  -I & I \\
\end{array}%
\right).
\end{equation*}
The hessian matrix of $w(x)$ is
$$ D^2w(x)=2 ke^{-k|x-y_0|^2}I-4
k^2e^{-k|x-y_0|^2}(x-y_0)\otimes(x-y_0).$$Its  eigenvalues are $2
ke^{-k|x-y_0|^2}$ with multiplicity $N-1$ and $2
ke^{-k|x-y_0|^2}(1-2k|x-y_0|^2)$ with multiplicity 1. In the
annulus $B_1\setminus \overline{B}_2$ we have $2
ke^{-k|x-y_0|^2}(1-2k|x-y_0|^2)\leq
2ke^{-k|x-y_0|^2}\left(1-k\frac{R^2}{2}\right)<0$ since
$k>\frac{2}{R^2}$.

Using the fact that $u$ and $v$ are respectively sub and
supersolution and the properties of the operator $F$ we have
\begin{equation*}\begin{split}g(y_j)&\leq F(y_j,v(y_j),-D_y\phi,Y_j)\\
&\leq F(y_j,v(y_j),-D_y\phi,Y_j+\frac{\sigma}{2} D^2w(y_j))
+\frac{\sigma}{2}\mathcal{M}_{a,A}^+( D^2w(y_j))\\& \leq
F(x_j,v(y_j),-D_y\phi,X_j-\frac{\sigma}{2}
D^2w(x_j))+\omega_T(o_j)+\frac{\sigma}{2} \mathcal{M}_{a,A}^+(
D^2w(y_j))\\ & \leq
F(x_j,u(x_j),D_x\phi,X_j)+\omega_T(o_j)+\frac{\sigma}{2}
\mathcal{M}_{a,A}^+( D^2w(x_j))+\frac{\sigma}{2}
\mathcal{M}_{a,A}^+(
D^2w(y_j))\\&+b\frac{\sigma}{2}|Dw(x_j)|+b\frac{\sigma}{2}|Dw(y_j)|+c|u(x_j)-v(y_j)|\\&
\leq g(x_j)+\omega_T(o_j)+\frac{\sigma}{2} \mathcal{M}_{a,A}^+(
D^2w(x_j))+\frac{\sigma}{2} \mathcal{M}_{a,A}^+(
D^2w(y_j))+b\frac{\sigma}{2}|Dw(x_j)|\\&+b\frac{\sigma}{2}|Dw(y_j)|+c|u(x_j)-v(y_j)|,
\end{split}\end{equation*}
where $o_j=j|x_j-y_j|^2+|x_j-y_j|(|D_y\phi|+1)\rightarrow0$ as
$j\rightarrow +\infty$. Then
\begin{equation*}\begin{split}g(y_j)&\leq g(x_j)+
A(N-1)\sigma ke^{-k|x_j-y_0|^2}+a\sigma
ke^{-k|x_j-y_0|^2}(1-2k|x_j-y_0|^2)\\&+A(N-1)\sigma
ke^{-k|y_j-y_0|^2}+a\sigma
ke^{-k|y_j-y_0|^2}(1-2k|y_j-y_0|^2)\\&+k\sigma
e^{-k|x_j-y_0|^2}b|x_j-y_0| +k\sigma e^{-k|y_j-y_0|^2}b|y_j-y_0|
+c|u(x_j)-v(y_j)|+\omega_T(o_j).\end{split}\end{equation*} Passing
to the limit as $j\rightarrow+\infty$ we get
\begin{equation*}\begin{split}2\sigma e^{-k|\zs-y_0|^2}\{-2ak^2|\zs-y_0|^2+[A(N-1)+a+b
|\zs-y_0|]k\}+c|u(\zs)-v(\zs)|\geq0.\end{split}\end{equation*}

Using \eqref{stcompstimau-v} and the fact that
$\frac{R}{2}<|\zs-y_0|\leq R$, we have
\begin{equation*}\begin{split}&0\leq 2\sigma e^{-k|\zs-y_0|^2}\{-2ak^2|\zs-y_0|^2+[A(N-1)+a+b
|\zs-y_0|]k\}+c|u(\zs)-v(\zs)|\\&\leq 2\sigma
e^{-k|\zs-y_0|^2}\left\{-ak^2\frac{R^2}{2}+[A(N-1)+a+b
R]k\right\}+\sigma c( e^{-k|\zs-y_0|^2}-e^{-kR^2})\\ &\leq \sigma
e^{-k|\zs-y_0|^2}\{-ak^2R^2+2[A(N-1)+a+b
R]k+c\}.\end{split}\end{equation*} If we fix $k>2/R^2$ so large
that $$-ak^2R^2+2[A(N-1)+a+b R]k+c<0,$$we  obtain a contradiction,
then $u<v$ on $\Oms$.\finedim
\begin{rem}{\em In Theorem \ref{stcompneu} the domain $\Om$ may be
unbounded. In that case, in the proof of the theorem it suffices
to maximize  $\psi(x,y)$ on the compact set
$(\overline{B}(y_0,2R)\cap\Oms)^2$, instead of the whole
$\Oms$.}\end{rem}A consequence of Theorem \ref{stcompneu} are the
following strong maximum and minimum principles.
\begin{cor}\label{strongmaxp}Assume the hypothesis of Theorem \ref{stcompneu}.
If $f(x,0)\leq 0$ for any $x\in\partial\Om$ and $v\in LSC(\Oms)$
is a non-negative viscosity supersolution of
\begin{equation}\label{stromaxpsis}
\begin{cases}
 F(x,v,Dv,D^2v)= 0  & \text{in} \quad\Om \\
 B(x,v,Dv)= 0 & \text{on} \quad\partial\Om, \\
 \end{cases}
 \end{equation} then either $v\equiv 0$ or $v>0$ on $\Oms$. If $f(x,0)\geq
 0$ for any $x\in\partial\Om$ and $u\in USC(\Oms)$
is a non-positive viscosity subsolution of \eqref{stromaxpsis}
then either $u\equiv 0$ or $u<0$ on $\Oms$.
\end{cor}
\dim If $f(x,0)\leq 0$ for any $x\in\partial\Om$ then $u\equiv 0$
is a subsolution of \eqref{stromaxpsis}. The thesis follows
applying Theorem \ref{stcompneu}. \finedim
\section{The Maximum Principle and the principal eigenvalues}
We say that  $F$ with boundary condition \eqref{e2} satisfies the
maximum principle, if whenever $u\in USC(\Oms)$ is a viscosity
subsolution of
\begin{equation*}
\begin{cases}
 F(x,u,Du,D^2u)= 0  & \text{in} \quad\Om \\
 B(x,u,Du)= 0 & \text{on} \quad\partial\Om, \\
 \end{cases}
 \end{equation*}then $u\leq 0$ on $\Oms$. We first prove that the
 maximum principle holds if  $F$ is proper, i.e., if $r\rightarrow F(x,r,p,M)$ is
non-decreasing. Observe that we do not require the stronger
condition $F(x,r,p,X) -\sigma r$ non-decreasing in $r$ for some
$\sigma>0$, in which case the comparison principle holds (see
\cite{cil} Theorem 7.5) and implies the maximum principle if
$u\equiv 0$ is a supersolution.

Successively, we show that the operator $F-\lam I$ with boundary
condition \eqref{e2} satisfies the maximum principle for any
$\lam<\lams$. To prove that the two results do not coincide, we
construct a class of operators which are not proper but that have
positive principal eigenvalue $\lams$, hence for which the maximum
principle holds.
\subsection{The case F proper}
\begin{thm}\label{macpcless0}Assume
that (F2), (F3), (f1) and (f3) hold, that $r\rightarrow
F(x,r,p,M)$ is non-decreasing on $\R$ for all $(x,p,M)\in
\Oms\times\R^N\times S(N)$, $F(x,0,0,0)\geq 0$ for all $x\in\Om$,
$f(x,0)\geq0$ for all $x\in\partial\Om$ and
\begin{equation}\label{t1e1ishii}\max_{x\in\partial\Om}
f(x,r)\vee\max_{x\in\Oms}F(x,r,0,0)>0 \text{ for any
}r>0.\end{equation}If $u\in USC(\overline{\Om})$ is a viscosity
subsolution of
\begin{equation}\label{t1e1}
\begin{cases}
 F(x,u,Du,D^2u)= 0  & \text{in} \quad\Om \\
 B(x,u,Du)= 0 & \text{on} \quad\partial\Om, \\
 \end{cases}
 \end{equation}
then $u\leq 0$ on $\Oms$.
\end{thm}
\dim Let $u$ be a subsolution of \eqref{t1e1}. First let us
suppose $u\equiv k=$const. By definition of subsolution and
Proposition \ref{ishii}
$$F(x,k,0,0)\leq 0, \text{ for any }x\in\Om$$ and $$B(x,k,0)=f(x,k)\leq0\text{ for any }x\in\partial\Om.$$
Then the hypothesis \eqref{t1e1ishii} implies $k\leq0$.

Now we assume that $u$ is not a constant. We argue by
contradiction; suppose that $\max_{\Oms}u=u(x_0)>0$, for some
$x_0\in\Oms$. Define $\widetilde{u}(x):=u(x)-u(x_0)$. Since
$r\rightarrow F(x,r,p,M)$ and $r\rightarrow f(x,r)$ are
non-decreasing, $\widetilde{u}$ is a non-positive subsolution of
\eqref{t1e1}. The properties $F(x,0,0,0)\geq 0$ and $f(x,0)\geq0$
imply that $v\equiv 0$ is a supersolution of \eqref{t1e1}. Then it
follows from Theorem \ref{stcompneu} that either $u\equiv u(x_0)$
or $u<u(x_0)$ on $\overline{\Om}$. In both cases we get a
contradiction.

\finedim
\begin{rem}{\em Under the assumptions of Theorem
\ref{macpcless0}, but now with $F(x,0,0,0)\leq 0$ for all
$x\in\Om$, $f(x,0)\leq0$ for all $x\in\partial\Om$ and
$\min_{x\in\partial\Om} f(x,r)\vee\min_{x\in\Oms}F(x,r,0,0)<0 $
for any $r<0$, we can prove the minimum principle, i.e., if $u\in
LSC(\Oms)$ is a viscosity supersolution of \eqref{t1e1} then
$u\geq 0$ on $\Oms$.}\end{rem}
\begin{rem}\emph{If $F$ does not depend on $r$ and $f\equiv 0$ a counterexample to the validity of the maximum
principle is given by the positive constants.}\end{rem}

\subsection{The Maximum Principle for $\lam<\lams$}

We set
\begin{equation*}\begin{split}\overline{E}:=\{&\lam\in
\R\;|\;\exists\,v>0 \text{ on $\Oms$ bounded viscosity
supersolution of }\\&
F(x,v,Dv,D^2v)=\lam v \text{ in } \Om,\, B(x,v,Dv) = 0\text{ on }\partial \Om \},\\
\end{split}\end{equation*}
\begin{equation*}\begin{split}\underline{E}:=\{&\lam\in
\R\;|\;\exists\,u<0 \text{ on $\Oms$ bounded  viscosity
subsolution of }\\&
F(x,u,Du,D^2u)=\lam u \text{ in } \Om,\, B(x,u,Du)= 0\text{ on }\partial \Om \}.\\
\end{split}\end{equation*}
The set $\overline{E}$ is not empty, indeed the function
$v(x)=e^{-|f(\cdot,1)|_\infty d(x)}$ satisfies
\begin{equation*}\begin{split}F(x,v,Dv,D^2v)-\lam v&\geq e^{-|f(\cdot,1)|_\infty d(x)}\big\{-\mathcal{M}_{a,A}^+
\big(|f(\cdot,1)|_\infty^2 Dd(x)\otimes Dd(x)\\&
-|f(\cdot,1)|_\infty D^2d(x)\big)-b|f(\cdot,1)|_\infty
-c-\lam\big\}\geq 0,\end{split}\end{equation*}in $\Om$, for $\lam$
small enough, and
\begin{equation*}\begin{split}B(x,v,Dv)&=f(x,1)+|f(x,1)|_\infty \geq
0,\end{split}\end{equation*}on $\partial\Om$. As a consequence
$\lams=\sup\overline{E}$ is well defined. Similarly we can prove
that $\underline{E}$ is not empty. We shall show that $\lams$ and
$\lamso$ are finite.

We want to remark that since in the sequel we will assume (f2),
which implies $f(x,0)=0$ for any $x\in\partial\Om$, by Corollary
\ref{strongmaxp} any non-negative supersolution (resp.,
non-positive subsolution) of $F(x,v,Dv,D^2v)=0$ in $\Om$,
$B(x,v,Dv)=0$ on $\partial\Om$ which is non-zero will be positive
(resp., negative) in all $\Oms$.
\begin{thm}\label{uequivv} Assume
that (F1)-(F3), (f1) and (f2) hold. Let $u\in USC(\Oms)$ and $v\in
LSC(\Oms)$ be respectively sub and supersolution of
\begin{equation*}
\begin{cases}
 F(x,u,Du,D^2u)= 0 & \text{in} \quad\Om \\
 B(x,u,Du)= 0 & \text{on} \quad\partial\Om. \\
 \end{cases}
 \end{equation*}
If $v$ is bounded, $v>0$ on $\Oms$ and $u(x_0)>0$ for some
$x_0\in\Oms$, then there exists $t>0$ such that $v\equiv tu$. The
same conclusion holds if $u$ is bounded, $u<0$ on $\Oms$ and
$v(x_0)<0$.
\end{thm}
\dim Suppose that $v>0$ on $\Oms$ and $u(x_0)>0$. We prove the
theorem trough a typical argument which is used in \cite{bnv} for
the linear case and Dirichlet boundary condition. Set $w_t=u-tv$.
If $t$ is large enough $w_t<0$ on $\Oms$. We define
$$\tau=\inf\{t\,|\,w_t<0 \text{ on } \Oms\}.$$ Clearly $w_{\tau}\leq 0$. If
$\max_{\Oms}w_{\tau}=m<0$, then for any $x\in\Oms$
$$w_{\tau-\epsilon}(x)=u(x)-(\tau-\epsilon)v(x)\leq m+\epsilon |v|_\infty<0,$$ for $\epsilon$
small enough. This contradicts the definition of $\tau$. Then
$w_{\tau}$ vanishes somewhere on $\Oms$ and $\tau>0$ since
$u(x_0)>0$. In conclusion $u\leq \tau v$ and $u(x)=\tau v(x)$ for
some $x\in\Oms$. Since $\tau v$ is again a supersolution, by
Theorem \ref{stcompneu} we have $u\equiv \tau v$.

If the inequalities satisfied by $u$ and $v$ are reversed, that is
$u<0$ and $v(x_0)<0$, we consider the function $w_t=tu-v$ and use
the same argument.\finedim
\begin{thm}[Maximum Principle for $\lam<\lams$]\label{maxpneum} Assume
that (F1)-(F3), (f1) and (f2) hold and $\lam<\lams$. Let $u\in
USC(\Oms)$ be a viscosity subsolution of
\begin{equation}\label{maxpeq}
\begin{cases}
 F(x,u,Du,D^2u)=\lam u & \text{in} \quad\Om \\
 B(x,u,Du)= 0 & \text{on} \quad\partial\Om, \\
 \end{cases}
 \end{equation}
then $u\leq 0$ on $\Oms$.
\end{thm}
\dim Let $\tau\in ]\lam,\lams[$, then by definition there exists
 $v>0$ on $\Oms$ bounded viscosity supersolution of
 \begin{equation*}
 \begin{cases}
 F(x,v,Dv,D^2v)=\tau v   & \text{in} \quad\Om \\
 B(x,v,Dv) =0 & \text{on} \quad\partial\Om. \\
 \end{cases}
 \end{equation*}
 Then $v$  satisfies
  \begin{equation}\label{th3.1ishii}
 \begin{cases}
 F(x,v,Dv,D^2v)-\lam v \geq (\tau-\lam)v> 0  & \text{in} \quad\Om \\
 B(x,v,Dv) \geq 0 & \text{on} \quad\partial\Om, \\
 \end{cases}
 \end{equation}in the viscosity sense.
 Suppose by contradiction that $u(x_0)>0$ for some $x_0\in\Oms$.
 Applying Theorem \ref{uequivv} to the operator $F-\lam I$, there
 exists $t>0$ such that $u\equiv t v.$ Then $u$ is positive on $\Oms$ and by homogeneity satisfies \eqref{th3.1ishii} in the viscosity sense. Since in addition $u$ is a viscosity subsolution of \eqref{maxpeq}, using Lemma 7.3 of \cite{iy} we get  $$(\tau-\lam)u\leq 0\quad\text{in }\Om,$$
 which is a contradiction.
\finedim
\begin{rem}{\em Similarly, we can prove the minimum principle for $\lam<\underline{\lam}$, i.e.,  if
$u\in LSC(\Oms)$ is a viscosity supersolution of \eqref{maxpeq}
and $\lam<\underline{\lam}$ then $u\geq 0$ on $\Oms$.}
\end{rem}
\begin{cor}Under the assumptions of Theorem \ref{maxpneum}, the quantities $\lams$ and $\underline{\lam}$ are
finite.
\end{cor}
\dim By Theorem \ref{maxpneum} it suffices to find $\lam\in\R$ and
a function $w$ which is a positive subsolution of
\begin{equation*}
 \begin{cases}
 F(x,w,Dw,D^2w)=\lam w  & \text{in} \quad\Om \\
 B(x,w,Dw) = 0 & \text{on} \quad\partial\Om. \\
 \end{cases}
 \end{equation*}For  \begin{equation*}\begin{split}\lam\geq -\mathcal{M}_{a,A}^-
 \left(|f(\cdot,1)|_\infty^2 Dd(x)\otimes Dd(x)+|f(\cdot,1)|_\infty D^2d(x)\right)+b|f(\cdot,1)|_\infty +c,\end{split}\end{equation*}
 a subsolution is  $w(x)=e^{|f(\cdot,1)|_\infty d(x)}$.
 \finedim
 \subsection{An example}
 We want to show  there exist some operators which are not proper but
whose first eigenvalue $\lams$ is positive.

For simplicity, let us suppose  that $F$ is independent of the
gradient variable and that $\Om$ is the ball of center 0 and
radius R. We assume in addition that for all
$(x,X)\in\Oms\times\emph{S(N)}$ and any $r>0$
\begin{equation}\label{ccambiasgn}F(x,r,X)\geq
-\mathcal{M}_{a,A}^+(X)+c_0(x)r,\end{equation} for some functions
$c_0(x)$. The Isaacs operator \eqref{isaacs} satisfies
\eqref{ccambiasgn} if
$$c_{\al,\beta}(x)\geq c_0(x)\quad \text{for all }x\in\Oms\text{ and }(\al,\beta)\in
\mathcal{A}\times\mathcal{B}.$$ In this case the operator is
proper if $c_0(x)\geq 0$. Since we are interested in non-proper
$F$, we are looking for functions  $c_0(x)$ in \eqref{ccambiasgn}
that may be negative somewhere. We suppose that
\begin{equation*}
\begin{cases}
c_0(x)>0 &\text{if } R-\epsilon<|x|\leq R\\
c_0(x)\geq \beta_1 & \text{if }\rho<|x|\leq R-\epsilon\\
c_0(x)\geq -\beta_2 & \text{if }|x|\leq \rho,\\
 \end{cases}
 \end{equation*}
 where $0<\rho<R$, $\epsilon> 0$ is small enough and
$\beta_1,\, \beta_2>0$. Remark that in the ball of radius $\rho$,
$c_0(x)$ may assume negative values. To prove that $\lams>0$ it
suffices to find $v>0$ bounded supersolution of
\begin{equation*}
\begin{cases}
 -\mathcal{M}_{a,A}^+(D^2v)+c_0(x)v=\lam v & \text{in} \quad\Om \\
 B(x,v,Dv)= 0 & \text{on} \quad\partial\Om, \\
 \end{cases}
 \end{equation*}for some $\lam>0$. Assume $f(x,r)\geq 0$ for any $x\in\partial\Om$ and $r\geq0$, then, as shown in \cite{p}, such supersolution
  $v$ exists if $\beta_1$ and $\beta_2$ satisfy the following inequality for
 some $k>0$
 \begin{equation*}\beta_2<\frac{ke^{-k\rho}a\left(k+\frac{N-1}{\rho}\right)}
{k\frac{R-\rho}{4}+k\frac{2NAR-(N-1)a(R+\rho)}{\beta_1R(R-\rho)}+1-e^{-k\rho}}.
\end{equation*}
As observed in \cite{p}, from the last relation we can see that
choosing $k=\frac{1}{\rho}$ the term on the right-hand side goes
to $+\infty$ as $\rho\rightarrow 0^+$, that is, if the set where
$c_0(x)$ is negative becomes smaller then the values of $c_0(x)$
in this set can be very negative. On the contrary, for any value
of $k$, if $\rho\rightarrow R^-$ then $\beta_2$ goes to 0.
Finally, for any $k$, if $\beta_1\rightarrow 0^+$ then again
$\beta_2$ goes to 0. So there is a sort of balance between
$\beta_1$ and $\beta_2$. In \cite{p} we present an example to
explain this behavior. For operators which satisfy
\eqref{ccambiasgn}, the property $\lams>0$ can be proved in any
$C^2$ domain, under similar assumptions on $c_0(x)$, see \cite{p}.
\section{Lipschitz regularity}
In this section we shall prove that viscosity solutions are
Lipschitz continuous on $\Oms$. We want to mention the works of Barles and Da Lio \cite{bdl} and Milakis and Silvestre \cite{ms} about H\"{o}lder estimates of viscosity solutions of fully nonlinear elliptic equations associated to Neumann type boundary conditions.
\begin{thm}\label{regolarita}Assume that (F1), (F2),
(F4), (f1) and (f2) hold. Let $g$ be a bounded function and $u\in
C(\Oms)$ be a viscosity solution of
\begin{equation*}
\begin{cases}
 F(x,u,Du,D^2u)=  g(x) & \text{in} \quad\Om \\
 B(x,u,Du)= 0 & \text{on} \quad\partial\Om, \\
 \end{cases}
 \end{equation*}
 then there exists $C_0>0$ such that
\begin{equation}\label{stimau-v1}|u(x)-u(y)|\leq C_0|x-y|\quad \forall
x,y\in\Oms,\end{equation}
 where $C_0$ depends on
 $N,\,a,\,A,\,b,\,c,\,C_1$, $\Om$ and $|f(\cdot,u(\cdot))|_\infty$.
 \end{thm}
 \dim
 We follow the proof of Proposition III.1 of \cite{il}, that we modify taking test functions which depend on the distance function and that are suitable
 for the Neumann boundary conditions.

We set
$$\Phi(x)=MK|x|-M(K|x|)^2,$$and
$$\varphi (x,y)=e^{-L(d(x)+d(y))}\Phi(x-y),$$ where
$L$ is a fixed number greater than $\frac{2}{3r}$ with $r$ the
radius in the condition ($\Om2$) and $K$ and $M$ are two positive
constants to be chosen later. If $K|x|\leq \frac{1}{4}$, then
\begin{equation}\label{phimagg}\Phi(x)\geq
\frac{3}{4}MK|x|.\end{equation} We define
$$\Delta _K:=\left\{(x,y)\in \R^N\times\R^N:\, |x-y|\leq\frac{1}{4K}\right\}.$$ We
fix $M$ such that
\begin{equation}\label{M}\max_{\Oms^{\,2}}|u(x)-u(y)|\leq
e^{-2Ld_0}\frac{M}{8},\end{equation} where
$d_0=\max_{x\in\Oms}d(x),$ and we claim that taking $\delta$ small
enough and $K$ large enough, one has
\begin{equation}\label{u-vindelta}
\delta( u(x)-u(y))-\varphi(x,y)\leq 0\quad\text{for }(x,y)\in
\Delta_K\cap\Oms^2.
\end{equation}
In this case \eqref{stimau-v1} is proven. To show
\eqref{u-vindelta} we suppose by contradiction that for some
$(\xs,\ys)\in\Delta_K\cap \Oms^2$
\begin{equation}\label{u-vcontr}\delta u(\xs)-\delta u(\ys)-\varphi(\xs,\ys)=\max_{\Delta_K\cap \Oms\,^2}(\delta u(x)-\delta u(y)-\varphi(x,y))>0.\end{equation}
Observe that $\delta u$ is again a solution since both $F$ and $B$
are positively homogeneous. Here we have dropped the dependence of
$\xs,\,\ys$ on $K$ and $\delta$ for simplicity of notations.

Clearly $\xs\neq \ys$. Moreover the point $(\xs,\ys)$ belongs to
$\text{int}(\Delta_K)\cap\Oms^2$. Indeed, if $|x-y|=\frac{1}{4K}$,
by \eqref{M} and \eqref{phimagg} for $\delta\leq 1$ we have
\begin{equation*}\begin{split}
\delta u(x)-\delta u(y)&\leq |u(x)-u(y)|\leq e^{-2Ld_0}
\frac{M}{8}\leq
e^{-L(d(x)+d(y))}\frac{1}{2}MK|x-y|\leq\varphi(x,y).\end{split}\end{equation*}
Since $\xs\neq \ys$ we can compute the derivatives of $\varphi$ in
$(\xs,\ys)$ obtaining
\begin{equation*}\begin{split}D_x\varphi(\xs,\ys)=&-Le^{-L(d(\xs)+d(\ys))}MK|\xs-\ys|(1-K|\xs-\ys|)Dd(\xs)\\&+
e^{-L(d(\xs)+d(\ys))}MK(1-2K|\xs-\ys|)\frac{(\xs-\ys)}{|\xs-\ys|},\end{split}\end{equation*}
\begin{equation*}\begin{split}D_y\varphi(\xs,\ys)=&-Le^{-L(d(\xs)+d(\ys))}MK|\xs-\ys|(1-K|\xs-\ys|)Dd(\ys)\\&-e^{-L(d(\xs)+d(\ys))}MK
(1-2K|\xs-\ys|)\frac{(\xs-\ys)}{|\xs-\ys|}.\end{split}\end{equation*}
 Observe that for $K\geq \frac{L}{4}$
\begin{equation}\label{dxphistima}
|D_x\varphi(\xs,\ys)|,|D_y\varphi(\xs,\ys)|\leq 2MK.\end{equation}
Using \eqref{sferaest}, if $\xs\in\partial \Om$ we have
\begin{equation}\label{condbordo}\begin{split}
B(\xs,\delta u(\xs),D_x\varphi(\xs,\ys))&=f(\xs,\delta
u(\xs))+Le^{-Ld(\ys)}MK|\xs-\ys|(1-K|\xs-\ys|)
\\&+e^{-Ld(\ys)}MK(1-2K|\xs-\ys|)\langle
\overrightarrow{n}(\xs),\frac{(\xs-\ys)}{|\xs-\ys|}\rangle \\&
\geq
\frac{1}{2}e^{-Ld(\ys)}MK|\xs-\ys|\left(\frac{3}{2}L-\frac{1}{r}\right)-\delta
|f(\cdot,u(\cdot))|_\infty>0,
\end{split}\end{equation}
since $\xs\neq\ys$, $L>\frac{2}{3r}$, for $\delta$ small enough.
Similarly, if $\ys\in\partial\Om$ then
\begin{equation*}
B(\ys,\delta u(\ys),-D_y\varphi(\xs,\ys))\leq
\frac{1}{2}e^{-Ld(\xs)}MK|\xs-\ys|\left(-\frac{3}{2}L+\frac{1}{r}\right)+\delta
|f(\cdot,u(\cdot))|_\infty<0.\end{equation*} Then $\xs,\ys\in\Om$
and
$$F(\xs,\delta u(\xs),D_x\varphi(\xs,\ys),X)\leq \delta g(\xs),\quad\text{if }(D_x\varphi(\xs,\ys),X)\in \overline{J}^{2,+}\delta
u(\xs),$$ $$F(\ys,\delta u(\ys),-D_y\varphi(\xs,\ys),Y)\geq \delta
g(\ys)\quad\text{if }(-D_y\varphi(\xs,\ys),Y)\in
\overline{J}^{2,-}\delta u(\ys).$$ Since
$(\xs,\ys)\in\text{int}\Delta_K\cap\Oms\,^2$, it is a local
maximum point of $\delta u(x)-\delta u(y)-\varphi(x,y)$ in
$\Oms\,^2$. Then applying Theorem 3.2 in \cite{cil}, for every
$\epsilon>0$ there exist $X,Y\in \emph{S(N)}$ such that $
(D_x\varphi(\xs,\ys),X)\in \overline{J}\,^{2,+}\delta
u(\xs),\,(-D_y\varphi(\xs,\ys),Y)\in \overline{J}\,^{2,-}\delta
u(\ys)$ and
\begin{equation}\label{tm2ishii}\left(%
\begin{array}{cc}
  X & 0 \\
  0 & -Y \\
\end{array}%
\right)\leq D^2(\varphi(\xs,\ys))+\epsilon
(D^2(\varphi(\xs,\ys)))^2.
\end{equation}
Now we want to estimate the matrix on the right-hand side of the
last inequality.
\begin{equation*}\begin{split}D^2\varphi(\xs,\ys)&=\Phi(\xs-\ys)D^2(e^{-L(d(\xs)+d(\ys))})+D(e^{-L(d(\xs)+d(\ys))})\otimes
D(\Phi(\xs-\ys))\\&+D(\Phi(\xs-\ys))\otimes
D(e^{-L(d(\xs)+d(\ys))})+e^{-L(d(\xs)+d(\ys))}D^2(\Phi(\xs-\ys)).\end{split}\end{equation*}We
set $$A_1:=\Phi(\xs-\ys)D^2(e^{-L(d(\xs)+d(\ys))}),$$
$$A_2:=D(e^{-L(d(\xs)+d(\ys))})\otimes
D(\Phi(\xs-\ys))+D(\Phi(\xs-\ys))\otimes
D(e^{-L(d(\xs)+d(\ys))}),$$
$$A_3:=e^{-L(d(\xs)+d(\ys))}D^2(\Phi(\xs-\ys)).$$Observe that
\begin{equation}\label{a1}A_1\leq CK|\xs-\ys|\left(\begin{array}{cc}
  I & 0 \\
  0 & I \\
\end{array}%
\right).\end{equation}  Here and henceforth C denotes various
positive constants independent of $K$ and $\delta$.

For $A_2$ we have the following estimate
\begin{equation}\label{a2}A_2\leq
CK\left(%
\begin{array}{cc}
  I & 0 \\
  0 & I \\
\end{array}%
\right)
+CK\left(%
\begin{array}{cc}
  I & -I \\
  -I & I \\
\end{array}%
\right).\end{equation} Indeed for $\xi,\,\eta\in\R^N$ we compute
\begin{equation*}\begin{split}\langle
A_2(\xi,\eta),(\xi,\eta)\rangle&=2Le^{-L(d(\xs)+d(\ys))}\{\langle
Dd(\xs)\otimes D\Phi(\xs-\ys)(\eta-\xi),\xi\rangle\\&+\langle
Dd(\ys)\otimes D\Phi(\xs-\ys)(\eta-\xi),\eta\rangle\}\leq
CK(|\xi|+|\eta|)|\eta-\xi|\\&\leq
CK(|\xi|^2+|\eta|^2)+CK|\eta-\xi|^2.\end{split}\end{equation*} Now
we consider $A_3$. The matrix $D^2(\Phi(\xs-\ys))$ has the form
$$D^2(\Phi(\xs-\ys))=\left(%
\begin{array}{cc}
  D^2\Phi(\xs-\ys) & - D^2\Phi(\xs-\ys) \\
  - D^2\Phi(\xs-\ys) &  D^2\Phi(\xs-\ys) \\
\end{array}%
\right),$$and the Hessian matrix of $\Phi(x)$ is
\begin{equation}\label{hessianphi}D^2\Phi(x)=\frac{MK}{|x|}\left(I-\frac{x\otimes
x}{|x|^2}\right)-2MK^2I.\end{equation} If we choose
\begin{equation}\label{epsilon}\epsilon=\frac{|\xs-\ys|}{2MKe^{-L(d(\xs)+d(\ys))}},\end{equation}
then we have the following estimates
$$\epsilon A_1^2\leq
CK|\xs-\ys|^3I_{2N},\quad \epsilon A_2^2\leq CK|\xs-\ys|I_{2N},$$
\begin{equation}\label{aprodotti}\begin{split} \epsilon (A_1A_2+A_2A_1)\leq
CK|\xs-\ys|^2I_{2N},\end{split}\end{equation}
$$\epsilon (A_1A_3+A_3A_1)\leq
CK|\xs-\ys|I_{2N},\quad \epsilon (A_2A_3+A_3A_2)\leq CKI_{2N},$$
 where $I_{2N}:=\left(%
\begin{array}{cc}
  I & 0 \\
  0 & I \\
\end{array}%
\right)$. Then using \eqref{a1}, \eqref{a2},
\eqref{aprodotti} and observing that $$(D^2(\Phi(\xs-\ys)))^2=\left(%
\begin{array}{cc}
  2(D^2\Phi(\xs-\ys))^2 & - 2(D^2\Phi(\xs-\ys))^2 \\
  - 2(D^2\Phi(\xs-\ys))^2 &  2(D^2\Phi(\xs-\ys))^2 \\
\end{array}%
\right),$$from \eqref{tm2ishii} we can conclude that
$$\left(%
\begin{array}{cc}
  X & 0 \\
  0 & -Y \\
\end{array}%
\right)\leq O(K)\left(%
\begin{array}{cc}
  I & 0 \\
  0 & I \\
\end{array}%
\right)+\left(%
\begin{array}{cc}
  B & -B \\
  -B & B \\
\end{array}%
\right),$$ where
\begin{equation}\label{matriceB}B=CKI+e^{-L(d(\xs)+d(\ys))}\left[D^2\Phi(\xs-\ys)+\frac{|\xs-\ys|}{MK}
(D^2\Phi(\xs-\ys))^2\right].\end{equation}The last inequality can
be rewritten as follows
$$\left(%
\begin{array}{cc}
  \widetilde{X} & 0 \\
  0 & -\widetilde{Y} \\
\end{array}%
\right)\leq\left(%
\begin{array}{cc}
  B & -B \\
  -B & B \\
\end{array}%
\right),$$ with $\widetilde{X}=X-O(K)I$ and
$\widetilde{Y}=Y+O(K)I.$

Now we want to get a good estimate for
tr($\widetilde{X}-\widetilde{Y}$), as in \cite{il}. For that aim
let
$$0\leq P:=\frac{(\xs-\ys)\otimes (\xs-\ys)}{|\xs-\ys|^2}\leq I.$$
Since $\widetilde{X}-\widetilde{Y}\leq 0$ and
$\widetilde{X}-\widetilde{Y}\leq 4B,$ we have
$$\text{tr}(\widetilde{X}-\widetilde{Y})\leq \text{tr}(P(\widetilde{X}-\widetilde{Y}))\leq 4 \text{tr}(PB).$$ We
have to compute tr($PB$). From \eqref{hessianphi}, observing that
the matrix $(1/|x|^2)x\otimes x$ is idempotent, i.e.,
$[(1/|x|^2)x\otimes x]^2=(1/|x|^2)x\otimes x$, we compute
$$(D^2\Phi(x))^2=\frac{M^2K^2}{|x|^2}(1-4K|x|)\left(I-\frac{x\otimes
x}{|x|^2}\right)+4M^2K^4I.$$ Then, since $\text{tr}P=1$ and
$4K|\xs-\ys|\leq1$, we have
\begin{equation*}\begin{split}\text{tr}(PB)&=CK+e^{-L(d(\xs)+d(\ys))}(-2MK^2+4MK^3|\xs-\ys|)\\&
\leq CK-e^{-L(d(\xs)+d(\ys))}MK^2<0, \end{split}\end{equation*}for
large $K$. This gives
$$|\text{tr}(\widetilde{X}-\widetilde{Y})|=-\text{tr}(\widetilde{X}-\widetilde{Y})\geq
4e^{-L(d(\xs)+d(\ys))}MK^2-4CK\geq CK^2,$$
 for large $K$. Since $\|B\|\leq \frac{CK}{|\xs-\ys|},$ we have
\begin{equation*}\begin{split}\|B\|^{\frac{1}{2}}|\text{tr}(\widetilde{X}-\widetilde{Y})|^{\frac{1}{2}}&\leq
\left(\frac{CK}{|\xs-\ys|}\right)^{\frac{1}{2}}|\text{tr}(\widetilde{X}-\widetilde{Y})|^{\frac{1}{2}}
\leq
\frac{C}{K^{\frac{1}{2}}|\xs-\ys|^{\frac{1}{2}}}|\text{tr}(\widetilde{X}-\widetilde{Y})|.
\end{split}\end{equation*}The Lemma III.I in \cite{il} ensures the
existence of a universal constant $C$ depending only on $N$ such
that $$\|\widetilde{X}\|, \|\widetilde{Y}\|\leq
C\{|\text{tr}(\widetilde{X}-\widetilde{Y})|+\|B\|^{\frac{1}{2}}|\text{tr}(\widetilde{X}-\widetilde{Y})|^{\frac{1}{2}}\}.$$
Thanks to the above estimates we can conclude that
\begin{equation}\label{normaxy}\|\widetilde{X}\|,\,\|\widetilde{Y}\|\leq
C|\text{tr}(\widetilde{X}-\widetilde{Y})|\left(1+\frac{1}{K^{\frac{1}{2}}|\xs-\ys|^{\frac{1}{2}}}\right).\end{equation}

Now, using assumptions (F2) and (F4) concerning $F$, the
definition of $\widetilde{X}$ and $\widetilde{Y}$ and the fact
that $\delta u$ is sub and supersolution we
compute\begin{equation*}\begin{split}\delta g(\ys)&\leq
F(\ys,\delta u(\ys),-D_y\varphi,Y)\leq F(\ys,\delta
u(\ys),-D_y\varphi,\widetilde{Y})+O(K)\\ & \leq F(\ys,\delta
u(\xs),D_x\varphi,\widetilde{X}) + c\delta|u(\xs)-u(\ys)|
+b|D_x\varphi+D_y\varphi|\\&+a\text{tr}(\widetilde{X}-\widetilde{Y})+O(K)\\&
\leq F(\xs,\delta
u(\xs),D_x\varphi,\widetilde{X})+2c\delta|u(\xs)|
+2b|D_x\varphi|+C_1|\xs-\ys|^\frac{1}{2}\|\widetilde{X}\|\\&+
c\delta|u(\xs)-u(\ys)|
+b|D_x\varphi+D_y\varphi|+a\text{tr}(\widetilde{X}-\widetilde{Y})+O(K)\\&\leq
\delta g(\xs)+2c\delta|u(\xs)|
+2b|D_x\varphi|+C_1|\xs-\ys|^\frac{1}{2}\|\widetilde{X}\|+
c\delta|u(\xs)-u(\ys)|\\&
+b|D_x\varphi+D_y\varphi|+a\text{tr}(\widetilde{X}-\widetilde{Y})+O(K).
\end{split}\end{equation*} From this inequalities, using
\eqref{dxphistima} and \eqref{normaxy}, we get
\begin{equation}\label{ultimalemm}\begin{split}
&\delta g(\ys)-\delta
g(\xs)-2c\delta|u(\xs)|-c\delta|u(\xs)-u(\ys)|\\&\leq
O(K)+C|\text{tr}(\widetilde{X}-\widetilde{Y})|(|\xs-\ys|^{\frac{1}{2}}+K^{-\frac{1}{2}})
+a\text{tr}(\widetilde{X}-\widetilde{Y})
\\&=a\text{tr}(\widetilde{X}-\widetilde{Y})+o(|\text{tr}(\widetilde{X}-\widetilde{Y})|),
\end{split}\end{equation}as $K\rightarrow+\infty$.
Since $g$ and $u$ are bounded, the first member in
\eqref{ultimalemm} is bounded from below by the quantity
$-2|g|_\infty-4c|u|_\infty$ which is independent of $\delta$. But
the last term in \eqref{ultimalemm} goes to $-\infty$ as
$K\rightarrow+\infty$, hence taking $K$ so large that
$$a\text{tr}(\widetilde{X}-\widetilde{Y})+o(|\text{tr}(\widetilde{X}-\widetilde{Y})|)<-2|g|_\infty-4c|u|_\infty,$$
and then $\delta$ so small that the last member in
\eqref{condbordo} is positive, we obtain a contradiction and this
concludes the proof.\finedim
\begin{rem}{\em The regularity theorem can be shown also for solutions of
the Neumann problem for the operator
\begin{equation*}\sup_{\al\in
\mathcal{A}}\inf_{\beta\in \mathcal{B}}\{
-\text{tr}(A_{\al,\beta}(x) D^2u)+b_{\al,\beta}(x)\cdot
Du+c_{\al,\beta}(x) u-g_{\al,\beta}(x)\},\end{equation*}if the
functions $g_{\al,\beta}$ are bounded uniformly in $\al$ and
$\beta$.}
\end{rem}
Since the Lipschitz estimate  depends only on the bounds of the
solution of $g$ and on the structural constants, an immediate
consequence of the previous theorem is the following compactness
criterion that will be useful in the next sections.
 \begin{cor}\label{corcomp}Assume the same hypothesis of Theorem \ref{regolarita}. Suppose that $(g_n)_{n}$ is a
sequence of continuous and uniformly bounded functions and
$(u_n)_n$ is a sequence of uniformly bounded viscosity solutions
of
\begin{equation*}
\begin{cases}
 F(x,u_n,Du_n,D^2u_n)=  g_n(x) & \text{in} \quad\Om \\
 B(x,u_n,Du_n)= 0 & \text{on} \quad\partial\Om. \\
 \end{cases}
 \end{equation*}Then the sequence $(u_n)_n$ is relatively compact in
 $C(\Oms)$.
 \end{cor}

\section{Existence results}
 This section is devoted to the
problem of the existence of a solution of
\begin{equation}\label{existence}
\begin{cases}
 F(x,u,Du,D^2u)=\lam u+  g(x) & \text{in} \quad\Om \\
 B(x,u,Du)= 0 & \text{on} \quad\partial\Om.\\
 \end{cases}
 \end{equation}
 Using the well known result which guarantees that \eqref{existence} with $\lam=0$ is uniquely solvable if $F$ satisfies
\begin{itemize}
\item[(F5)] There exists $\sigma>0$ such that for any $(x,p,X)\in\Oms\times\R^N\times\emph{S(N)}$
the function $r\rightarrow F(x,r,p,X)- \sigma r$ is non-decreasing
on $\R$,\end{itemize} see \cite{cil} Theorem 7.5,
 we will prove the
 existence of a positive solution of \eqref{existence} when $g$ is
 non-negative and $\lam<\lams$, without requiring (F5). The solution is unique if $g>0$.
 Then  we will show the existence of a positive principal
 eigenfunction corresponding to $\lams$, that is a solution of
 \eqref{existence} when $g\equiv0$ and $\lam=\lams$. For the last
 two results we will follow the proof given in \cite{bd} for the
 analogous theorems with the Dirichlet boundary condition.

 Symmetrical results can be obtained for the eigenvalue
 $\underline{\lam}$.

 Finally, we will prove that the Neumann problem  \eqref{existence} is solvable for any right-hand side if $\lam<\min\{\lams,\underline{\lam}\}$.

The following is a well known result, see \cite{cil} Theorem 7.5.
 \begin{thm}\label{princonfc<0}Suppose that (F2), (F3), (F5), (f1) and (f3) hold and that $g$ is
continuous on $\Oms$. If $u\in USC(\Oms)$ and $v\in LSC(\Oms)$ are
respectively sub and supersolution of
\begin{equation}\label{princonfl<ls2}
\begin{cases}
 F(x,u,Du,D^2u)=  g(x) & \text{in} \quad\Om \\
 B(x,u,Du)= 0 & \text{on} \quad\partial\Om, \\
 \end{cases}
 \end{equation}
then $u\leq v$ on $\Oms.$ Moreover \eqref{princonfl<ls2} has a
unique viscosity solution.
\end{thm}
\begin{thm}\label{comparisonl<ls}
Assume that (F1)-(F3), (f1) and (f2) hold. Suppose $h\geq0$,
$g\leq h$ and $g(x)<0$ if $h(x)=0$. Let $u\in USC(\Oms)$ be a
viscosity subsolution of
\begin{equation*}
\begin{cases}
 F(x,u,Du,D^2u)=\lam u + g(x) & \text{in} \quad\Om \\
 B(x,u,Du)= 0 & \text{on} \quad\partial\Om, \\
 \end{cases}
 \end{equation*}and let $v\in LSC(\Oms)$ be a bounded
positive viscosity supersolution of
\begin{equation}\label{compla<lams}
\begin{cases}
 F(x,v,Dv,D^2v)=\lam v+  h(x) & \text{in} \quad\Om \\
 B(x,v,Dv)= 0 & \text{on} \quad\partial\Om. \\
 \end{cases}
 \end{equation} Then $u\leq v$ on $\Oms.$
 \end{thm}
 \begin{rem}{\em The existence of such a $v$ implies
 $\lam\leq\lams$.}\end{rem}
 \begin{rem}\label{complamso}{\em Similarly, we can prove the comparison result between $u$ and $v$ if $u$ is negative and
 bounded,
 $g\leq0$, $g\leq h$ and $h(x)>0$ if $g(x)=0$.}\end{rem}
 \dim
 Suppose by contradiction that $\max_{\Oms}(u-v)=u(\xs)-v(\xs)>0$ for some $\xs\in\Oms$. Set $w_t=u-tv$. If $t$ is large enough $w_t<0$ on $\Oms$. We
define
$$\tau=\inf\{t\,|\,w_t<0 \text{ on } \Oms\}.$$As in the proof of Theorem \ref{uequivv}, $w_{\tau}\leq 0$ and vanishes in some point, i.e.,
$u\leq \tau v$ and $u(x)=\tau v(x)$ for some $x\in\Oms$. Moreover,
since $u(\xs)>v(\xs)$ we know that $\tau >1$, which implies that
$h\leq\tau h$, being $h$  non-negative. Then $\tau v$ is still a
supersolution of \eqref{compla<lams} and $u\equiv \tau v$ by
Theorem \ref{stcompneu}. Hence, applying Lemma 7.3 of \cite{iy} we get
$$  \tau h\leq g,$$ which contradicts the
assumptions on $g$ and $h$. \finedim
\begin{thm}\label{esistl<ls} Suppose that (F1)-(F4), (f1)-(f3) hold, that $\lam<\lams$, $g\geq 0$, $g\not\equiv0$ and $g$ is
continuous on $\Oms$, then there exists a positive viscosity
solution of \eqref{existence}. The positive solution is unique if
$g>0$.
\end{thm}
\dim The condition (F2) implies that $r\rightarrow F(x,r,p,X)+cr$
is non-decreasing. Hence the operator $F+(2c+|\lam|)I$ satisfies
(F5) with $\sigma=c+|\lam|$, so that by Theorem \ref{princonfc<0}
the sequence $(u_n)_n$ defined by $u_1=0$ and $u_{n+1}$ as the
solution of
\begin{equation*}
\begin{cases}
 F(x,u_{n+1},Du_{n+1},D^2u_{n+1})+(2c+|\lam|)u_{n+1}=g+(2c+|\lam|+\lam)u_n& \text{in} \quad\Om \\
 B(x,u_{n+1}, Du_{n+1})= 0 & \text{on} \quad\partial\Om, \\
 \end{cases}
 \end{equation*}
is well defined. By the comparison Theorems \ref{princonfc<0} and
\ref{stcompneu}, since $g\geq 0$ and $g\not\equiv 0$ the sequence
is positive and increasing.

We use the argument of Theorem 7 of \cite{bd} to prove that
$(u_n)_n$ is also bounded. Suppose that it is not, then dividing
by $|u_{n+1}|_\infty$ and defining $v_n:=\frac{u_n}{|u_n|_\infty}$
one gets that $v_{n+1}$ is a solution of
\begin{equation*}
\begin{cases}
 F(x,v_{n+1},Dv_{n+1},D^2v_{n+1})+(2c+|\lam|)v_{n+1}\\\quad =\frac{g}{|u_{n+1}|_\infty}+(2c+|\lam|+\lam)\frac{u_n}{|u_{n+1}|_\infty}&
 \text{in} \quad\Om \\
 B(x,v_{n+1}, Dv_{n+1})= 0 & \text{on} \quad\partial\Om. \\
 \end{cases}
 \end{equation*}
 By Corollary \ref{corcomp},  $(v_n)_n$ converges along a subsequence to a positive function $v$ which
 satisfies
 \begin{equation*}
\begin{cases}
 F(x,v,Dv,D^2v)-\lam v=(2c+|\lam|+\lam)(k-1)v\leq 0& \text{in} \quad\Om \\
 B(x,v, Dv)= 0 & \text{on} \quad\partial\Om, \\
 \end{cases}
 \end{equation*}where
 $k:=\limsup_{n\rightarrow+\infty}\frac{|u_n|_\infty}{|u_{n+1}|_\infty}\leq
 1$. This contradicts the maximum principle, Theo\-rem
 \ref{maxpneum}.
 Then $(u_n)_n$ is bounded and letting $n$ go to infinity, by the compactness result, the sequence
converges uniformly to a function $u$ which is a solution.
Moreover the solution is positive on $\Oms$ by Corollary
\ref{strongmaxp}.

The uniqueness of the positive solution follows from Theorem
\ref{comparisonl<ls}. \finedim

\begin{thm}[Existence of principal eigenfunctions]\label{esistautof}
 Suppose that (F1)-(F4), (f1)-(f3) hold. Then there exists $\phi>0$ on $\Oms$ viscosity solution of
\begin{equation}\label{firsteig}
\begin{cases}
 F(x,\phi,D\phi,D^2\phi)=\lams\phi & \text{in} \quad\Om \\
 B(x,\phi,D\phi)= 0 & \text{on} \quad\partial\Om. \\
 \end{cases}
 \end{equation}Moreover $\phi$ is Lipschitz continuous on $\Oms$.
\end{thm}
\dim Let $\lam_n$ be an increasing sequence which converges to
$\lams$. Let $u_n$ be a positive solution of
\begin{equation*}
\begin{cases}
 F(x,u_n,Du_n,D^2u_n)=\lam_n u_n+1 & \text{in} \quad\Om \\
 B(x,u_n,Du_n)= 0 & \text{on} \quad\partial\Om. \\
 \end{cases}
 \end{equation*}
 By  Theorem \ref{esistl<ls} the sequence $(u_n)_n$ is well
 defined. Following the argument of the proof of Theorem 8 of
 \cite{bd} we can prove that it is unbounded, otherwise one would
 contradict the definition of $\lams$.
 Then, up to subsequence, $|u_n|_\infty\rightarrow+\infty$ as
 $n\rightarrow+\infty$ and defining $v_n:=\frac{u_n}{|u_n|_\infty}$
 one gets that $v_n$ satisfies
 \begin{equation*}
\begin{cases}
 F(x,v_n,Dv_n,D^2v_n)=\lam_n v_n+ \frac{1}{|u_n|_\infty} & \text{in} \quad\Om \\
 B(x,v_n,Dv_n)= 0 & \text{on} \quad\partial\Om. \\
 \end{cases}
 \end{equation*}
 Then, by Corollary \ref{corcomp}, we can extract a subsequence converging to a
 function $\phi$ with $|\phi|_\infty=1$ which is positive on $\Oms$ by Corollary \ref{strongmaxp} and is the desired solution.
 By  Theorem \ref{regolarita} the solution is also Lipschitz continuous
 on
 $\Oms$. \finedim
 \begin{rem}\label{thmslsot}{\em With the same arguments used in the proofs of Theorems  \ref{esistl<ls} and \ref{esistautof}
 one can prove: the existence of a negative viscosity solution of \eqref{existence}, for
$\lam<\underline{\lam}$ and $g\leq 0$, $g\not\equiv0$, which is
unique if $g<0$ by Remark \ref{complamso}; the existence of a
negative Lipschitz
 principal eigenfunction corresponding to $\underline{\lam}$,
  i.e., a solution of
\begin{equation}\label{firsteig2}
\begin{cases}
 F(x,\phi,D\phi,D^2\phi)=\underline{\lam}\phi & \text{in} \quad\Om \\
 B(x,\phi,D\phi)= 0 & \text{on} \quad\partial\Om. \\
 \end{cases}
 \end{equation}}
 \end{rem}
 \begin{thm}Suppose that (F1)-(F4), (f1)-(f3) hold. Suppose that $\lam<\min\{\lams,\underline{\lam}\}$ and
$g$ is continuous on $\Oms$, then there exists a viscosity
solution of \eqref{existence}.
\end{thm}
 \dim If $g\equiv 0$, by the maximum and minimum principles the only solution is $u\equiv 0$. Let us suppose $g\not\equiv0$.
 Since $\lam<\min\{\lams,\underline{\lam}\}$ by Theorem
 \ref{esistl<ls} and Remark \ref{thmslsot} there exist
 $v_0\in C(\Oms)$ positive viscosity solution of
 \eqref{existence} with right-hand side $|g|_\infty$ and $u_0\in C(\Oms)$ negative viscosity solution
 of \eqref{existence} with right-hand side $-|g|_\infty$.

 Let $(u_n)_n$ be the sequence defined in the proof of Theorem \ref{esistl<ls} with $u_1=u_0$. By
 comparison Theorem \ref{princonfc<0} we have $u_0=u_1\leq u_2\leq
 ...\leq v_0$. Hence, by the compactness Corollary \ref{corcomp}
 the sequence converges to a continuous function which is the
 desired solution.\finedim
 \begin{rem}{\em The existence results can be shown also for the operator
 \begin{equation*}
 \sup_{\al\in
\mathcal{A}}\inf_{\beta\in \mathcal{B}}\{
-\text{tr}(A_{\al,\beta}(x) D^2u)+b_{\al,\beta}(x)\cdot
Du+c_{\al,\beta}(x) u-g_{\al,\beta}(x)\},\end{equation*}if the
functions $g_{\al,\beta}$ are continuous uniformly in $\al$ and
$\beta$. In particular, in that case, if $\lams$ and $\lamso$ are
positive there exists a viscosity solution of
\begin{equation*}
\begin{cases}
 \sup_{\al\in\mathcal{A}}\inf_{\beta\in \mathcal{B}}\{-\text{tr}(A_{\al,\beta}(x) D^2u)+b_{\al,\beta}(x)\cdot
Du+c_{\al,\beta}(x) u-g_{\al,\beta}(x)\}=0\,\,  \text{in } \Om \\
 B(x,u,Du)= 0\quad\quad\quad\quad\quad\quad\quad\quad\quad\quad\quad\quad\quad\quad\quad\quad\quad\quad
 \quad\quad\quad\quad\quad\quad\quad\,\text{on } \partial\Om. \\
 \end{cases}
 \end{equation*}}
\end{rem}

\section{Properties of the principal eigenvalues}
In this section we establish some of the basic properties of the
principal eigenvalues. We denote by $\phi^+$ a positive
eigenfunction corresponding to $\lams$ and by $\phi^-$ a negative
eigenfunction corresponding  to $\lamso$. Throughout this section
we assume (F1)-(F4) and (f1)-(f3).

The next result states that the principal eigenfunctions are
simple, in the sense that they are equal up to a multiplicative
constant.
 \begin{prop}\label{simpleeig}If  $u\in USC(\Oms)$ is a viscosity
subsolution of
 \begin{equation}\label{simple2}
\begin{cases}
 F(x,u,Du,D^2u)=\lams u & \text{in} \quad\Om \\
 B(x,u,Du)= 0 & \text{on} \quad\partial\Om, \\
 \end{cases}
 \end{equation}and $u(x_0)>0$ for some $x_0\in\Oms$ then there exists $t>0$ such that $u\equiv
 t\phi^+$. If $u\in LSC(\Oms)$ is a viscosity supersolution of \eqref{simple2} with $\lams$ replaced by
 $\underline{\lam}$ and $u(x_0)<0$, then there exists $t>0$ such that $u\equiv t\phi^-$.

 Assume in addition \begin{equation}\label{simple3}-F(x,-r,-p,-X)\leq F(x,r,p,X)\text{ for any }
 (x,r,p,X)\in\Oms\times\R\times\R^N\times\emph{S(N)}\end{equation}
  and \begin{equation}\label{simple4}-f(x,-r)\leq f(x,r)\text{ for any }(x,r)\in\partial\Om\times\R.\end{equation}
  If $u\in C(\Oms)$, $u\not\equiv 0$, is a viscosity subsolution of \eqref{simple2} then there exists $t\in\R$ such that $u\equiv
 t\phi^+$. If $u\in C(\Oms)$, $u\not\equiv 0$ is a viscosity solution of \eqref{simple2} with $\lams$ replaced by $\underline{\lam}$,
 there exists $t\in\R$ such that $u\equiv t\phi^-$.
\end{prop}
\dim If $u$ is a subsolution (resp., supersolution) of
\eqref{simple2} (resp., of \eqref{simple2} with $\underline{\lam}$
instead of $\lams$) and $u(x_0)>0$ (resp., $u(x_0)<0$), then by
Theorem \ref{uequivv} we have $u\equiv
 t\phi^+$ (resp., $u\equiv
 t\phi^-$) for some $t>0$.

Now assume \eqref{simple3}-\eqref{simple4} and let $u\not\equiv0$
be a subsolution of \eqref{simple2}. If $u$ is positive somewhere
we are in the previous case. If $u$ is negative on $\Oms$ then the
function $w:=-u$ is a positive continuous supersolution of
 \begin{equation*}
\begin{cases}
 F(x,w,Dw,D^2w)-\lams w\geq -F(x,-w,-Dw,-D^2w)+\lams(-w)\geq 0 & \text{in} \quad\Om \\
 B(x,w,Dw)\geq-B(x,-w,-Dw)\geq 0 & \text{on} \quad\partial\Om. \\
 \end{cases}
 \end{equation*}  Hence, again
from Theorem \ref{uequivv} it follows that $u\equiv t\phi^+$, for
some $t<0$.

Finally, let $u\not\equiv0$ be a solution of \eqref{simple2} with
$\underline{\lam}$ instead of $\lams$. Remark that conditions
\eqref{simple3}-\eqref{simple4} imply $\underline{\lam}\leq\lams$.
If $\underline{\lam}<\lams$, then by the maximum principle,
Theorem \ref{maxpneum}, $u<0$ on $\Oms$ and we are in the first
case. If $\underline{\lam}=\lams$, by the simplicity of $\lams$
just proved, $u\equiv t \phi^-$ for some $t<0$. \finedim
\begin{rem}{\em If $F$ and $f$ satisfy
\begin{equation*}-F(x,-r,-p,-X)\geq F(x,r,p,X)\text{ for any }
 (x,r,p,X)\in\Oms\times\R\times\R^N\times\emph{S(N)}\end{equation*}
  and \begin{equation*}-f(x,-r)\geq f(x,r)\text{ for any
  }(x,r)\in\partial\Om\times\R,\end{equation*} then, applying
  Proposition \ref{simpleeig} to the operator
  $G(x,r,p,X)=-F(x,-r,-p,-X)$ with
  $B(x,r,p)=\widetilde{f}(x,r)+\langle
  p,\overrightarrow{n}(x)\rangle$, where
  $\widetilde{f}(x,r)=-f(x,-r)$, we get again simplicity of
  principal eigenvalues.}\end{rem}
\begin{rem}{\em Convex and 1-homogeneous operators satisfy the assumption \eqref{simple3}.}\end{rem}
\begin{prop}\label{nootherpositive}$\lams$ (resp., $\underline{\lam}$)
is the only eigenvalue corresponding to a positive (resp.,
negative) eigenfunction.
\end{prop}
\dim Let $u$ be a positive eigenfunction corresponding to $\mu$. By the definition of $\lams$, we have $\mu\leq\lams$. If $\mu<\lams$, we must have $u\leq 0$ by Theorem \ref{maxpneum}, which is a contradiction. Thus $\mu=\lam$.\finedim The following
proposition states that the principal eigenvalues are isolated.
 \begin{prop}There exists $\epsilon>0$ such that the
problem
 \begin{equation}\label{isolati}
\begin{cases}
 F(x,u,Du,D^2u)=\lam u & \text{in} \quad\Om \\
 B(x,u,Du)= 0 & \text{on} \quad\partial\Om, \\
 \end{cases}
 \end{equation}has no solutions $u\not\equiv 0$, for
 $\lam\in(-\infty,\max\{\underline{\lam},\lams\}+\epsilon)\setminus\{\lams,\underline{\lam}\}$.
 \end{prop}
 \dim
 We may suppose
without loss of generality that $\lams\leq \underline{\lam}$.
If $\lam<\lams\leq \underline{\lam}$ then it follows from  the
maximum and minimum principles that $u\equiv 0$ is the only
 solution of \eqref{isolati}.

 If $\lam<\underline{\lam}$ and $u\not\equiv 0$
 is a solution of \eqref{isolati}, by the minimum principle we
 have
 $u>0$ on $\Oms$. Then Proposition \ref{nootherpositive} implies $\lam=\lams$.

 Finally suppose that there exists a sequence $\lam_n\downarrow\underline{\lam}$ such that the problem \eqref{isolati} with $\lam=\lam_n$ has
 a solution $\phi_n\not\equiv0$. We can assume that
 $|\phi_n|_\infty=1$ for any $n$. Then by the compactness criterion,
 Corollary \ref{corcomp}, the sequence $(\phi_n)_n$ converges uniformly on
 $\Oms$ to a function $\phi\not\equiv0$ which is a solution of
 \eqref{isolati} with $\lam=\underline{\lam}$. By Proposition \ref{nootherpositive} the functions $\phi_n$ change sign in $\Om$ while by
 Proposition
 \ref{simpleeig} and Theorem \ref{stcompneu}  either
 $\phi>0$ or $\phi<0$ on $\Oms$.  This contradicts the uniform convergence of
 $(\phi_n)_n$ to $\phi$.\finedim
 We want to conclude this section with the following comparison, suggested by Hitoshi Ishii \cite{icom},
  between  $\lams=\lams_N$ and $\lams_D$ respectively the principal eigenvalues
 corresponding to the Neumann and the Dirichlet problems.
 \begin{prop}$\lams_N<\lams_D$.\end{prop}
 \dim Let $v$ and $w$ be respectively the eigenfunctions
 corresponding to $\lams_N$ and $\lams_D$. That is
$$F(x,v,Dv,D^2v)= \lams_N v \text{ in } \Om,\quad B(x,v,Dv)= 0  \text{ on }
\partial\Om,\quad v>0\text{ on } \Oms,$$
$$F(x,w,Dw,D^2w)= \lams_D w \text{ in } \Om,\quad w= 0  \text{ on }
\partial\Om,\quad w>0\text{ in } \Om.$$
Since $f(x,0)=0$, we see that $w$ satisfies
$$F(x,w,Dw,D^2w)= \lams_D w \text{ in } \Om,\quad B(x,w,Dw)\leq 0  \text{ on }
\partial\Om.$$ Let us suppose $\lams_N\geq\lams_D$. Then
$$F(x,w,Dw,D^2w)\leq\lams_N w \text{ in } \Om,\quad B(x,w,Dw)\leq 0  \text{ on }
\partial\Om.$$Replacing $w$ by its constant multiple $tw$ with
$t>0$, we may assume that $w\leq v$ on $\Oms$ and $w(x_0)=v(x_0)$
for some $x_0\in\Om$. Note that $w(x)=0<v(x)$ for all
$x\in\partial\Om$. By Theorem \ref{puccistcomdir} we must have
$w\equiv v$ or $w<v$ on $\Oms$. This is a contradiction. \finedim

\section{The Pucci's operators}
In this section we want to show that the two principal eigenvalues
of the following operator
$$F(x,u,Du,D^2u)=-\mathcal{M}_{a,A}^+(D^2u)+b(x)\cdot
Du+c(x)u,$$ with the pure Neumann boundary condition may be
different. Suppose $b\in C^{0,1}(\Oms)$, $c\in C^{0,\beta}(\Oms)$
for some $\beta>0$ and $\Om$ of class $C^{2,\beta}$.

If $ c(x)\equiv c_0$ is constant then it is easy to see that
$\lams=\underline{\lam}=c_0$ and by Proposition \ref{simpleeig}
the only eigenfunctions are the constants. Nevertheless, if $c(x)$
is not constant the two principal eigenvalues never coincide,
unless $\mathcal{M}_{a,A}^+$ is the Laplacian. To prove this we
need the following lemma, whose proof is given for the sake of
completness.
\begin{lem}\label{lapllemma}Suppose that $\Om$ is a $C^{2,\beta}$ domain,
 $b\in C^{0,\beta}(\Oms)$  and $c\in
C^{0,\beta}(\Oms)$, for some $0<\beta\leq 1$. Then the viscosity solutions of
\begin{equation}\label{classeig}
\begin{cases}
 -\Delta u+b(x)\cdot Du+c(x)u=0& \text{in} \quad\Om \\
 \langle Du, \overrightarrow{n}(x)\rangle= 0 & \text{on}\quad \partial\Om, \\
 \end{cases}
 \end{equation}are in $C^2(\Oms)$.
 \end{lem}
 \dim
 Consider the problem
\begin{equation}\label{lemma8.1}
\begin{cases}
 -\Delta v+b(x)\cdot Dv+v
  =f(x) & \text{in } \Om \\
 \langle Dv, \overrightarrow{n}(x)\rangle= 0 & \text{on } \partial\Om, \\
 \end{cases}
 \end{equation}
 where $f(x)=(1-c(x))u(x)$. By Theorem \ref{regolarita}, $u$ is Lipschitz continuous on $\Oms$ and then the function $f$ is H\"{o}lder continuous on $\Oms$. Moreover, it
is clear that $u$ is a solution of \eqref{lemma8.1}. The classical theory says that \eqref{lemma8.1} has a solution
$v\in C^2(\Oms)$. By uniqueness of viscosity solutions of \eqref{lemma8.1}, we find that $u = v$.
\finedim
\begin{prop}Assume the hypothesis of Lemma \ref{lapllemma} and let $b\in C^{0,1}(\Oms)$.
If $A\neq a$ and $\lams=\underline{\lam}$ then $c(x)$ is constant.
\end{prop}
\dim Let $\phi$ be a positive eigenfunction of $\lams$, i.e.
\begin{equation}\label{proplams}
\begin{cases}
 -\mathcal{M}_{a,A}^+(D^2\phi)+b(x)\cdot D\phi+(c(x)-\lams)\phi= 0 & \text{in} \quad\Om \\
 \langle D\phi, \overrightarrow{n}(x)\rangle= 0 & \text{on} \quad\partial\Om, \\
 \end{cases}
 \end{equation}
and let $-\psi$ be a negative eigenfunction corresponding to
$\underline{\lam}$. Since
$\mathcal{M}_{a,A}^+(-D^2\psi)=-\mathcal{M}_{a,A}^-(D^2\psi)$,
$\psi$ satisfies
\begin{equation}
\begin{cases}\label{proplams2}
 -\mathcal{M}_{a,A}^-(D^2\psi)+b(x)\cdot D\psi+(c(x)-\underline{\lam})\psi=  0& \text{in} \quad\Om \\
 \langle D\psi, \overrightarrow{n}(x)\rangle= 0 & \text{on} \quad\partial\Om. \\
 \end{cases}
 \end{equation}
 If $\lams=\underline{\lam}$ then by
 Proposition \ref{simpleeig} $\psi=t\phi$ for some $t>0$. We can
 assume $\psi=\phi$. By summing the first equations in
 \eqref{proplams} and \eqref{proplams2}, we can see that $\phi$ is a positive viscosity solution of
\begin{equation*}
\begin{cases}
 -(A+a)\Delta \phi+2b(x)\cdot D\phi+2(c(x)-\lams)\phi= 0 & \text{in} \quad\Om \\
 \langle D\phi, \overrightarrow{n}(x)\rangle= 0 & \text{on} \quad\partial\Om. \\
 \end{cases}
 \end{equation*}Then by Lemma
 \ref{lapllemma}, $\phi\in C^2(\Oms)$. Subtracting the first equations in
 \eqref{proplams} and \eqref{proplams2}, we can see that $\phi$ is a classical solution of
 $$(A-a)\sum_{i=1}^N|e_i(x)|= 0\quad\text{in} \quad\Om ,$$ where $e_1(x),...,e_N(x)$ are the eigenvalues of $D^2\phi(x)$.  Since $A\neq a$, the last equation implies that $e_i(x)=0$ in $\Om$ for any $i=1...N$. In particular, taking into consideration the boundary condition, $\phi$ is a classical solution of
 \begin{equation*}
\begin{cases}
 \Delta \phi= 0 & \text{in} \quad\Om \\
 \langle D\phi, \overrightarrow{n}(x)\rangle= 0 & \text{on} \quad\partial\Om, \\
 \end{cases}
 \end{equation*}
 and then has to be constant. This implies
 $$c(x)-\lams=0\quad \text{in }\Om,$$i.e., $c\equiv\lams$ is constant.\finedim

\begin{ack}{\em The author wishes to thank the Professor Hitoshi Ishii  for the fruitful discussions about the topics of this
paper.}\end{ack}


\begin{thebibliography}{10}


\bibitem {a} {\sc H. Amann}, Fixed point equations and nonlinear
eigenvalue problems in ordered Banach spaces, {\em SIAM Rev.},
{\bf 18} (1976), 620-709.

\bibitem {an} {\sc A. Anane}, Simplicité et isolation de la
première valeur propre du p-Laplacien avec poids. (French)
[Simplicity and isolation of the first eigenvalue of the
p-Laplacian with weight] {\em C. R. Acad. Sci. Paris Sr I Math.},
{\bf 305} (1987), no. 16, 752-728.

\bibitem {bc} {\sc M. Bardi and I. Capuzzo Dolcetta}, Optimal
control and viscosity solutions of Hamilton-Jacobi-Bellmann
equations, Birkhauser, 1997.


\bibitem {b}
{\scshape G. Barles}, Nonlinear Neumann boundary conditions for
quasilinear degenerate elliptic equations and applications,  {\em
J. Differential Equations}, {\bf 154} (1999), 191-224.

\bibitem {bdl}
{\scshape G. Barles and F. Da Lio}, Local $C^{0,\al}$ estimates for viscosity solutions of
Neumann-type boundary value problems,  {\em
J. Differential Equations}, {\bf 225} (2006),  no. 1, 202-241.


\bibitem{bnv} {\sc H. Berestycki, L. Nirenberg and S.R.S. Varadhan},
The principal eigenvalue and maximum principle for second order
elliptic operators in general domain, {\em Comm. Pure Appl.
Math.}, {\bf 47} (1994), no. 1, 47-92.

\bibitem {bd2}
{\scshape I. Birindelli and F. Demengel}, Comparison principle and
Liouville type results for singular fully nonlinear operators,
\emph{Ann. Fac. Sci. Toulouse Math.} {\bf 13} (2004), no. 2,
261-287.

\bibitem {bd}
{\scshape I. Birindelli and F. Demengel}, Eigenvalue, maximum
principle and regularity for fully nonlinear homogeneous
operators, \emph{Comm. Pure Appl. Anal.}, {\bf 6} (2007), no. 2,
335-366.


\bibitem{beq} {\sc J. Busca, M. J. Esteban, A. Quaas}, Nonlinear
eigenvalues and bifurcation problems for Pucci's operators, {\em
Ann. Inst. H. Poincaré Anal. Non Linéaire}, {\bf 22} (2005), no.
2, 187-206.

\bibitem{cc} {\sc L. Caffarelli and X. Cabré} Fully nonlinear
equations Colloquium Publications 43, American Mathematical
Society, Providence, RI, 1995.

\bibitem{cil} {\sc M.C. Crandall, H. Ishii and P.L. Lions},
User's guide to viscosity solutions of second order partial
differential equations, {\em Bull. Amer. Math. Soc. (N.S.)}, {\bf
27} (1992), no. 1, 1-67.

\bibitem{fs} {\sc W. H. Fleming and H. M. Soner}, Controlled
Markov processes and viscosity solutions. Applications of
Mathematics (New York), {\bf 25}. Springer-Verlag, New York, 1993.

\bibitem{go} {\sc Y. Giga and M. Ohnuma} On strong comparison
principle for semicontinouos viscosity solutions of some nonlinear
elliptic equations, \emph{preprint}.

\bibitem{gt}{\sc D. Gilbarg and N.S. Trudinger}, Elliptic partial
differential equations of second order. Reprint of the 1998
edition. Classics in Mathematics. Springer-Verlag, Berlin, 2001.

\bibitem{i} {\sc H. Ishii},
Fully nonlinear oblique derivative problems for Nonlinear
Second-Order Elliptic PDE's, {\em Duke Math. J.}, {\bf 62} (1991),
no. 3, 633-661.

\bibitem{i2} {\sc H. Ishii}, Perron's method for Hamilton-Jacobi
Equations,  {\em Duke Math. J.}, {\bf 55} (1987), 369-384.

\bibitem{icom} {\sc H. Ishii}, Personal Communication.

\bibitem{il} {\sc H. Ishii and P.L. Lions},
Viscosity Solutions of Fully Nonlinear Second-Order Elliptic
Partial Differential Equations, {\em  J. Differential Equations},
{\bf 83} (1990), no. 1, 26-78.



\bibitem{iy} {\sc H. Ishii and Y. Yoshimura}, {\em
Demi-eigenvalues for uniformly elliptic Isaacs operators},
preprint.

\bibitem{l} {\sc P. L. Lions}, Bifurcation and optimal stochastic
control, {\em Nonlinear Anal.}, {\bf 7} (1983), no. 2, 177-207.


\bibitem{li2} {\sc P. Lindqvist}, On a nonlinear eigenvalue problem,
{\em Fall. School in Analysis} (Jyv\"{a}skyl\"{a}, 1994), 33-54.
Report, 68, {\em Univ. Jyv\"{a}skyl\"{a}}, Jyv\"{a}skyl\"{a},
1995.

\bibitem{ms} {\sc E. Milakis and L. E. Silvestre}, Regularity for fully nonlinear elliptic equations
with Neumann boundary data, {\em Comm. Partial Differential Equations}, {\bf 31}
(2006), no. 7-9, 1227-1252.

\bibitem{p} {\sc S. Patrizi}, The Neumann problem for singular fully nonlinear
operators, \emph{to appear in Journal de Mathématiques Pure et
Appliquées}.

\bibitem{pw} {\sc M.H. Protter and H.F. Weinberger}, Maximum
principles in differential equations. Prentice-Hall, Inc.,
Englewood Cliffs, N.J. 1967.

\bibitem{pw2} {\sc M.H. Protter and H.F. Weinberger}, On the
spectrum of general second order operators, {\em Bull. AMS}, {\bf
72} (1966), 251-255.

\bibitem{q} {\sc A. Quaas}, Existence of positive solutions to a
"semilinear" equation involving the Pucci's operators in a convex
domain, {\em Differential Integral Equations}, {\bf 17}, (2004),
no. 5-6, 481-494.

\bibitem{qs} {\sc A. Quaas and B. Sirakov}, On the principal
eigenvalues and the Dirichlet problem for fully nonlinear
operators, {\em C. R. Math. Acad. Sci. Paris}, {\bf 342}, no. 2,
115-118.

\end{thebibliography}
\end{document}